\documentclass[12pt]{article}
\usepackage{latexsym,amssymb,amsmath}
\begin{document}

\baselineskip=20pt

\newcommand{\rd}{\mbox{Rad}}
\newcommand{\kn}{\mbox{ker}}
\newcommand{\psp}{\vspace{0.4cm}}
\newcommand{\pse}{\vspace{0.2cm}}
\newcommand{\ptl}{\partial}
\newcommand{\dlt}{\delta}
\newcommand{\Dlt}{\Delta}
\newcommand{\sgm}{\sigma}
\newcommand{\al}{\alpha}
\newcommand{\be}{\beta}
\newcommand{\G}{\Gamma}
\newcommand{\gm}{\gamma}
\newcommand{\lmd}{\lambda}
\newcommand{\td}{\tilde}
\newcommand{\vf}{\varphi}
\newcommand{\ad}{\mbox{ad}}
\newcommand{\stl}{\stackrel}
\newcommand{\ol}{\overline}
\newcommand{\es}{\epsilon}
\newcommand{\ves}{\varepsilon}
\newcommand{\la}{\langle}
\newcommand{\ra}{\rangle}
\newcommand{\vt}{\vartheta}
\newcommand{\wt}{\mbox{wt}\:}
\newcommand{\sym}{\mbox{sym}}
\newcommand{\for}{\mbox{for}}

\begin{center}{\Large \bf Theta Series of Unimodular Lattices, Combinatorial }\end{center}
\begin{center}{\Large \bf Identities and Weighted Symmetric Polynomials}\end{center}
\vspace{0.2cm}

\begin{center}{\large Xiaoping Xu}\end{center}
\begin{center}{Department of Mathematics, The Hong Kong University of Science \& Technology}\end{center}
\begin{center}{Clear Water Bay, Kowloon, Hong Kong}\footnote{Research supported
 by the Direct Allocation Grant DAG99/00.SC25 from HKUST}\end{center}

\vspace{0.3cm}

\begin{center}{\Large \bf Abstract}\end{center}
\vspace{0.2cm}

{\small We find two  combinatorial identities on  the theta series of the root lattices of the finite-dimensional simple Lie algebras of type $D_{4n}$ and the cosets in their integral duals, in terms of the well-known Essenstein series $E_4(z)$ and Ramanujan series $\Dlt_{24}(z)$. Using these two identities, we determine the theta series of certain infinite families of postive definite even unimodular lattices obtained by gluing finite copies of the root lattices of the finite-dimensional simple Lie algebras of type $D_{2n}$. It turns out that these theta series are weighted symmetric polynomials of two fixed families of polynomials of $E_4(z)$ and $\Dlt_{24}(z)$}.
 
\section{Introduction}

Denote by $\Bbb{Z}$ the ring of integers. A unimodular (linear) lattice is a finite-rank free abelian group with a symmetric integer-valued $\Bbb{Z}$-bilinear form whose associated
symmetric matrix is of determinant $\pm 1$. Unimodular lattices are important objects in geometry of numbers. In sphere packings, many packings provided by unimodular lattices are of higher covering density (cf. [CS3]). The second cohomology group of a 4-manifold over $\Bbb{Z}$ modulo the torsion subgroup  forms a unimodular lattice with respect to its intersection matrix (cf. [F]).  The vertex operator superalgebra associated with a positive definite unimodular lattice enjoys the property of having a unique locally-finite irreducible module, that is, itself. The well-known example of  Leech lattice is directly related to 
important simple finite groups, the Conway groups (cf. [C]). The moonshine representation of the Monster group was constructed by orbifold construction through the vertex operator algebra associated with the Leech lattice and some of their twisted modules (cf. [B2], [FLM]). Classification of positive definite unimodular lattices has been done up to rank twenty-six (cf. [B1], [CS2]). In [X1] (also cf. [X2], [X3]), the author gave some explicit constructions of  infinite families of positive definite unimodular lattices.

 The arithmatic content of a positive definite unimodular lattice is given by its theta series, the generating function of the numbers of lattice points on spheres of integral square radius. A classical Hecke's theorem says that the theta series of a positive definite even unimodular lattice is a polynomial of the well-known  Essenstein series $E_4(z)$ and Ramanujan series $\Dlt_{24}(z)$ (cf. [H]). However, there is still a little known on how to write the theta series of a positive definite even  unimodular lattice as the polynomials of $E_4(z)$ and $\Dlt_{24}(z)$. In particular, there are no results on the theta series of infinite families of positive definite unimodular lattices.  
A natural question is what kind polynomials of $E_4(z)$ and $\Dlt_{24}(z)$ could be the theta series of positive definite even unimodular lattices. 
In this paper, we find two combinatorial identities on the theta series of the root lattices of the finite-dimensional simple Lie algebras of type $D_{4n}$ and the cosets in their integral duals, in terms of the Essenstein series $E_4(z)$ and Ramanujan series $\Dlt_{24}(z)$. Using these two identities, we determine the theta series of certain infinite families of postive definite even unimodular lattices with a sublattice of the same rank that is isomorphic to the direct sum of finite copies of the root lattices of the finite-dimensional simple Lie algebras of type $D_{2n}$. It turns out that these theta series are weighted symmetric polynomials in $\Dlt_{24}(z)$ and  two fixed families of polynomials of $E_4(z)$ and $\Dlt_{24}(z)$. Below, we shall give a more detailed technical introduction.

The fundamental arithematic functions used in counting lattice points are the following theta functions:
$$ \vt_2(z)=\sum_{m\in\Bbb{Z}}q^{(m+1/2)^2},\qquad\vt_3(z)=\sum_{m\in\Bbb{Z}}q^{m^2},\qquad\vt_4(z)=\sum_{m\in\Bbb{Z}}(-q)^{m^2},\eqno(1.1)$$
where
$$q=e^{\pi\sqrt{-1}z}\;\;\mbox{with}\;\; \mbox{Im}(z)>0.\eqno(1.2)$$
 An important identity among them is
$$\vt_2^4(z)+\vt_4^4(z)=\vt_3^4(z).\eqno(1.3)$$
The Essenstein $E_4(z)$ series is defined by
$$E_4(z)={1\over 2}(\vt_2^8(z)+\vt_3^8(z)+\vt_4^8(z))\eqno(1.4)$$
and the Ramanujan series is defined by
$$\Dlt_{24}(z)=\left({\vt_2(z)\vt_3(z)\vt_4(z)\over 2}\right)^8.\eqno(1.5)$$
As $q$-powers, 
$$E_4(z)=1+240\sum_{m=1}^{\infty}(\sum_{d|m}d^3)q^{2m}\eqno(1.6)$$
and 
$$\Dlt_{24}(z)=\sum_{m=1}^{\infty}\tau(m)q^{2m}\eqno(1.7)$$
whose coefficients $\tau(m)$ are called {\it Ramanujan numbers}. We refer to Table 4.9 in [CS3] for the first one hundred of coefficients of the $q$-powers in $E_4(z)$ and the reference [L] for first three hundred of the Ramanujan numbers.

Denote by $\Bbb{R}$ the field of real numbers. For a positive integer $n$, we denote by $\Bbb{R}^n$ the $n$-dimensional Eucleadean space with the inner product
$$\la\vec{\al},\vec{\be}\ra=\sum_{i=1}^n\al_i\be_i\qquad\for\;\;\vec{\al}=(\al_1,...,\al_n),\vec{\be}=(\be_1,...,\be_n)\in\Bbb{R}^n.\eqno(1.8)$$
A {\it positive definite integral lattice} $L$ {\it of rank} $n$ is an additive subgroup $L$ of $\Bbb{R}^n$ generated by a basis of $\Bbb{R}^n$ such that
$$\la\vec{\al},\vec{\be}\ra\in\Bbb{Z}\qquad\for\;\;\vec{\al},\vec{\be}\in L.\eqno(1.9)$$
The lattice $L$ is called {\it even} if
$$\la\vec{\al},\vec{\al}\ra\in 2\Bbb{Z}\qquad\for\;\;\vec{\al}\in L,\eqno(1.10)$$
and is unimodular if and only if
$$L=\{\vec{\al}\in \Bbb{R}^n\mid \la\vec{\al},\vec{\be}\ra\in\Bbb{Z}\;\for\;\vec{\be}\in L\}.\eqno(1.11)$$
For any coset $L'$ of $L$ in $\Bbb{R}^n$, we define the {\it theta series of} $L'$ by
$$\Theta_{L'}(z)=\sum_{\vec{\al}\in L'}q^{\la\vec{\al},\vec{\al}\ra}.\eqno(1.12)$$
Hecke [H] proved that $\Theta_L(z)$ is a polynomial of $E_4(z)$ and $\Dlt_{24}(z)$ for any positive definite even unimodular lattice $L$. Thus a fundamental question in positive definite even unimodular lattices is to determine the coefficients of the theta series of the lattices as polynomials of $E_4(z)$ and $\Dlt_{24}(z)$. One may ask the question in another way that what kind polynomials of $E_4(z)$ and $\Dlt_{24}(z)$ could be the theta series of positive definite even unimodular lattices. To this author's best knowledge, only a few of theta series of positive definite even unimodular lattices are known. 

 For a nonnegative integer $n$, we let
$$h_n(z)=\vt_2^{8n}(z)+\vt_3^{8n}(z)+\vt_4^{8n}(z),\eqno(1.13)$$
$$\rho_n(z)={\vt_3^{8(n+1)+4}(z)-\vt_2^{8(n+1)+4}(z)-\vt_4^{8(n+1)+4}(z)\over (\vt_2(z)\vt_3(z)\vt_4(z))^4}.\eqno(1.14)$$
By (1.3), the coefficients of the $q$-powers in $\rho_n(z)$ are positive integers.
Denote by $\Bbb{Z}_+$ the set of positive integers and by $\Bbb{N}$ the set of nonnegative integers.
In this paper, we shall prove that for any $n\in\Bbb{Z}_+$,  the following combinatorial identities hold:
$$h_n(z)=2E_4^n(z)+\sum_{i=1}^{[|n/3|]}{n\over i}\left(\!\!\!\begin{array}{c}n-i-1\\ 2i-1\end{array}\!\!\!\right)2^{8i}\Dlt_{24}^i(z)E_4^{n-3i}(z),\eqno(1.15)$$
$$\rho_n(z)=\sum_{i=0}^{[|n/3|]}{2n+3\over 2i+1}\left(\!\!\!\begin{array}{c}n-i\\ 2i \end{array}\!\!\!\right)2^{8i}\Dlt_{24}^i(z)E_4^{n-3i}(z).\eqno(1.16)$$
The above two identities can be viewed as higher-order analogues of the identity (1.3).

For a postive intger $n$, the type-D root lattice  $R_{D_n}$ is defined by
$$R_{D_n}=\{\al=(\al_1,\al_2,...,\al_n)\in\Bbb{Z}^n\mid \sum_{i=1}^n\al_i\in 2\Bbb{Z}\}\subset \Bbb{R}^n\eqno(1.17)$$
with the symmetric form (1.8). We define {\it weights} by
$$\wt\Dlt_{24}(z)=3,\;\;\wt h_n(z)=\wt \rho_n(z)=n\qquad \for\;\;n\in\Bbb{Z}_+.\eqno(1.18)$$
In this paper, we shall determine the theta series of certain infinite families of even unimodular lattices that have a sublattice of the same rank and isomorphic to the direct sum of finite copies of the lattices $R_{D_n}$ with various $n$. They are the following functions
\begin{eqnarray*}& &2^{-2\ell-1}[\sum_{1\leq j_1\leq j_2;j_1+j_2\leq \ell}\sym\{h_{m_1+\cdots +m_{2j_1}}(z)h_{m_{2j_1+1}+\cdots+m_{2(j_1+j_2)}}(z)h_{m_{2(j_1+j_2)+1}+\cdots+m_{2\ell+1}}(z)\}\\& &-\sum_{0\leq j_1\leq j_2\leq \ell-j_1-j_2-1}\sym\{h_{m_1+\cdots+m_{2j_1+1}}(z)h_{m_{2j_1+2}+\cdots+m_{2(j_1+j_2+1)}}(z)\\& &\times h_{m_{2(j_1+j_2)+3}+\cdots+m_{2\ell+1}}(z)\}
+\sum_{j=1}^{\ell}(3-4^{\ell-j})\sym\{h_{m_1+\cdots +m_{2j}}(z) h_{m_{2j+1}+\cdots +m_{2\ell+1}}(z)\}\\&&+\left(4-3\cdot 4^{\ell}+\frac{2}{3}\sum_{j_1+j_2\leq \ell-1}\left(\!\!\!\begin{array}{c}2\ell+1\\ 2j_1+1,2j_2+1\end{array}\!\!\!\right)\right)h_{m_1+\cdots+ m_{2\ell+1}}(z)]\hspace{3cm}(1.19)\end{eqnarray*}
for $\ell\in\Bbb{N}$ and $m_1,...,m_{2\ell+1}\in\Bbb{Z}_+$;
\begin{eqnarray*}& &2^{-2\ell}[\sum_{1\leq j_1\leq j_2\leq\ell-j_1-j_2}\sym\{h_{m_1+\cdots+m_{2j_1}+j_1}(z)h_{m_{2j_1+1}+\cdots+m_{2(j_1+j_2)}+j_2}(z)\\& &\times h_{m_{2(j_1+j_2)+1}+\cdots+m_{2\ell}+\ell-j_1-j_2}(z)\}
-2^8\Dlt_{24}(z)\sum_{1\leq j_1,j_2;j_1+j_2\leq\ell-2}\sym\{\rho_{m_1+\cdots+m_{2j_1+1}+j_1-1}(z)\\& &\times\rho_{m_{2j_1+2}+\cdots+m_{2(j_1+j_2+1)}+j_2-1}(z)h_{m_{2(j_1+j_2)+3}+\cdots+m_{2\ell}+\ell-j_1-j_2-1}(z)\}\hspace{4cm}\end{eqnarray*}
\begin{eqnarray*}& &+3\sum_{j=1}^{[|\ell/2|]}\sym\{h_{m_1+\cdots+m_{2j}+j}(z)
 h_{m_{2j+1}+\cdots+m_{2\ell}+\ell-j}(z)\}+2^8\Dlt_{24}(z)\\& &\times\sum_{j=0}^{[|(\ell-1)/2|]}(2^{2j}+2^{2(\ell-j-1)}-3)\sym\{\rho_{m_1+\cdots+m_{2j+1}+j-1}(z)\rho_{m_{2j+2}+\cdots+m_{2\ell}+\ell-j-2}(z)\}\\& &+\left(4-\frac{2}{3}\sum_{j_1+j_2\leq\ell}\left(\!\!\!\begin{array}{c}2\ell\\ 2j_1,2j_2\end{array}\!\!\!\right)\right)h_{m_1+\cdots+m_{2\ell}+\ell}(z)]\hspace{6cm}(1.20)\end{eqnarray*}
for $\ell\in\Bbb{Z}_+$ and $m_1,...,m_{2\ell}\in\Bbb{N}$;
\begin{eqnarray*}\hspace{1cm}& &{1\over 2}h_{m_1+m_2+m_3+m_4+2\es+1}(z)-32\Dlt_{24}(z)[\rho_{m_1+m_2+\es-1}(z)\rho_{m_3+m_4+\es-1}(z)\\& &+\rho_{m_1+m_3+\es-1}(z)\rho_{m_2+m_4+\es-1}(z)+\rho_{m_1+m_4+\es-1}(z)\rho_{m_2+m_3+\es-1}(z)]\hspace{2cm}(1.21)\end{eqnarray*}
for $m_1,m_2,m_3,m_4\in\Bbb{N}$ and $\es=0,1$.
 Here we have treated $\rho_{-1}(z)=0$ in (1.20), and (1.21), and ``$\sym\{\cdot\}$" means the symmetric polynomial with respect to the integrable variables $m_i$ and the expression in the braces is a representative term. For instance,
\begin{eqnarray*}& &\sym\{h_{m_1+m_2+1}(z)h_{m_3+m_4+1}(z)\}\\& =&h_{m_1+m_2+1}(z)h_{m_3+m_4+1}(z)+h_{m_1+m_3+1}(z)h_{m_2+m_4+1}(z)\\& &+h_{m_1+m_4+1}(z)h_{m_2+m_3+1}(z),\hspace{8.9cm}(1.22)\end{eqnarray*}
\begin{eqnarray*}& &\sym\{h_{m_1+m_2+1}(z)\rho_{m_3-1}(z)\rho_{m_4-1}(z)\}\\&=&h_{m_1+m_2+1}(z)\rho_{m_3-1}(z)\rho_{m_4-1}(z)+h_{m_1+m_3+1}(z)\rho_{m_2-1}(z)\rho_{m_4-1}(z)\\& &+h_{m_1+m_4+1}(z)\rho_{m_2-1}(z)\rho_{m_3-1}(z)+h_{m_3+m_4+1}(z)\rho_{m_1-1}(z)\rho_{m_2-1}(z)\\& &+h_{m_2+m_4+1}(z)\rho_{m_1-1}(z)\rho_{m_3-1}(z)+h_{m_2+m_3+1}(z)\rho_{m_1-1}(z)\rho_{m_4-1}(z).\hspace{2.3cm}(1.23)\end{eqnarray*}
These theta series  are weighted symmetric polynomials of the series $\{\Dlt_{24}(z),h_n(z),\rho_n(z)\mid n\in\Bbb{Z}_+\}$.  Thus we essentially determine  the theta series of these lattices as polynomials of the well-known Essenstein series $E_4(z)$ and Ramanujan series $\Dlt_{24}(z)$. We speculate that the theta series of the other infinite families of positive definite even unimodular lattices may relate to the invariants of the other finite groups.

The results in this paper could be useful in study of modular forms and partition functions of the conformal field theories related to positive definite even unimodular lattices. The covering densities of the sphere packings of these lattices can be calculated by our formulae.

In Section 2, we shall prove (1.15) and (1.16). In Section 3, the theta series of our concerned lattice will be determined.

\section{Proofs of the Combinatorial Identities}

In this section, we shall present the proofs of the combinatorial identities (1.15) and (1.16).

Set
$$a=\vt_2^4(z),\qquad b=\vt_4^4(z).\eqno(2.1)$$
By (1.3), we have
$$\vt_3^4(z)=a+b.\eqno(2.2)$$
Moreover, (1.13) and (1.14) become
$$h_n=a^{2n}+b^{2n}+(a+b)^{2n},\;\;\rho_n={(a+b)^{2n+3}-a^{2n+3}-b^{2n+3}\over ab(a+b)}.\eqno(2.3)$$
We allow $n=0$ in the above equations. 

For convenience, we let
$$E=E_4(z)={1\over 2}(a^2+b^2+(a+b)^2)=a^2+b^2+ab,\eqno(2.4)$$
$$\Dlt=2^8\Dlt_{24}(z)=(\vt_2(z)\vt_3(z)\vt_4(z))^8=a^2b^2(a+b)^2.\eqno(2.5)$$
Note that
$$h_0=3,\;\;\;h_1=2E.\eqno(2.6)$$
Moreover,
\begin{eqnarray*}\hspace{2cm}E^2&=&(a^2+b^2+ab)^2\\&=&a^4+b^4+(ab)^2+2a^2b^2+2a^2\cdot ab+2ab\cdot b^2\\&=&a^4+b^4+3a^2b^2+2ab(a^2+b^2).\hspace{6cm}(2.7)\end{eqnarray*}
Hence,
\begin{eqnarray*}\hspace{2cm}h_2&=&a^4+b^4+(a+b)^4\\&=&a^4+b^4+a^4+b^4+4ab(a^2+b^2)+6a^2b^2\\&=&2[a^4+b^4+3a^2b^2+2ab(a^2+b^2)]\\&=&2E^2.\hspace{10.9cm}(2.8)\end{eqnarray*}

{\bf Lemma 2.1}. {\it For} $3\leq n\in\Bbb{Z}$, {\it we have}:
$$h_n=2Eh_{n-1}-E^2h_{n-2}+\Dlt h_{n-3}.\eqno(2.9)$$

{\it Proof}. Note that
\begin{eqnarray*}2Eh_{n-1}&=&(a^2+b^2+(a+b)^2)(a^{2(n-1)}+b^{2(n-1)}+(a+b)^{2(n-1)})\hspace{5cm}\end{eqnarray*}
\begin{eqnarray*}&=&a^{2n}+b^{2n}+a^2b^2(a^{2(n-2)}+b^{2(n-2)})+(a^2+b^2)(a+b)^{2(n-1)}\\& &+(a+b)^2(a^{2(n-1)}+b^{2(n-1)})+(a+b)^{2n}\\&=&h_n+a^2b^2(a^{2(n-2)}+b^{2(n-2)})+(a^2+b^2)(a+b)^{2(n-1)}\\& &+(a+b)^2(a^{2(n-1)}+b^{2(n-1)}).\hspace{8.8cm}(2.10)\end{eqnarray*}
Moreover,
\begin{eqnarray*}\hspace{2cm}& &(a+b)^2(a^2+b^2)+a^2b^2\\&=&(a^2+b^2+2ab)(a^2+b^2)+a^2b^2\\&=&(a^2+b^2)^2+2ab(a^2+b^2)+a^2b^2\\&=&a^4+b^4+2a^2b^2+2ab(a^2+b^2)+a^2b^2\\&=&E^2.\hspace{11.4cm}(2.11)\end{eqnarray*}
Hence, we obtain
\begin{eqnarray*}& &E^2h_{n-2}\\&=&((a+b)^2(a^2+b^2)+a^2b^2)(a^{2(n-2)}+b^{2(n-2)}+(a+b)^{2(n-2)})\\&=&(a+b)^2[a^{2(n-1)}+b^{2(n-1)}+a^2b^2(a^{2(n-3)}+b^{2(n-3)})]+a^2b^2(a^{2(n-2)}+b^{2(n-2)})\\& &+(a^2+b^2)(a+b)^{2(n-1)}+a^2b^2(a+b)^{2(n-2)}
\\&=&(a+b)^2(a^{2(n-1)}+b^{2(n-1)})+\Dlt(a^{2(n-3)}+b^{2(n-3)})+a^2b^2(a^{2(n-2)}+b^{2(n-2)})\\& &+(a^2+b^2)(a+b)^{2(n-1)}+\Dlt(a+b)^{2(n-3)}\\&=&(a+b)^2(a^{2(n-1)}+b^{2(n-1)})+a^2b^2(a^{2(n-2)}+b^{2(n-2)})\\& &+(a^2+b^2)(a+b)^{2(n-1)}+\Dlt h_{n-3}.\hspace{8cm}(2.12)\end{eqnarray*}
Thus we have
$$2Eh_{n-1}-E^2h_{n-2}=h_n-\Dlt h_{n-3},\eqno(2.13)$$
which is equivalent to (2.9). $\qquad\Box$
\psp

By Lemma 2.1, we calculate
$$h_3=2E^3+3\Dlt,\;\;h_4=2E^4+8\Dlt E,\;\;h_5=2E^5+15\Dlt E^2,\eqno(2.14)$$
$$h_6=2E^6+24\Dlt E^3+3\Dlt^2,\;\;h_7=2E^7+35\Dlt E^4+14\Dlt^2 E,\eqno(2.15)$$
$$h_8=2E^8+48\Dlt E^5+40\Dlt^2E^2,\;\;h_9=2E^9+63\Dlt E^6+90\Dlt^2 E^3+3\Dlt^3,\eqno(2.16)$$
$$h_{10}=2E^{10}+80\Dlt E^7+175\Dlt^2E^4+20\Dlt^3 E.\eqno(2.17)$$
In fact, we had directly calculated (2.14)-(2.17), and then observed (2.9) from them. Analyzing the coefficients in (2.14)-(2.17), we speculated the following theorem.
\psp

{\bf Theorem 2.2}. {\it For} $n\in\Bbb{Z}_+$, {\it we have}
$$h_n=2E^n+\sum_{i=1}^{[|n/3|]}{n\over i}\left(\!\!\!\begin{array}{c}n-i-1\\ 2i-1 \end{array}\!\!\!\right)\Dlt^iE^{n-3i},\eqno(2.18)$$
{\it which is equivalent to (1.15) by (2.4) and (2.5)}.
\pse

{\it Proof}. We shall prove (2.18) by (2.9) and induction on $n$. First for $m,k\in\Bbb{Z}_+$ such that $k\leq m-1$, we have
\begin{eqnarray*}& &{2m\over k}\left(\!\!\!\begin{array}{c}m-k-1\\ 2k-1 \end{array}\!\!\!\right)-{m-1\over k}\left(\!\!\!\begin{array}{c}m-k-2\\ 2k-1 \end{array}\!\!\!\right)+{m-2\over k-1}\left(\!\!\!\begin{array}{c}m-k-2\\ 2k-3 \end{array}\!\!\!\right)\\&=&{1\over 2k(k-1)(2k-1)}\left(\!\!\!\begin{array}{c}m-k-2\\ 2k-3 \end{array}\!\!\!\right)[2m(m-k-1)(m+1-3k)\\& &-(m-1)(m+1-3k)(m-3k)+2k(m-2)(2k-1)]\\&=&
{1\over 2k(k-1)(2k-1)}\left(\!\!\!\begin{array}{c}m-k-2\\ 2k-3 \end{array}\!\!\!\right)[2m(m+1)(m-k-1)-6mk(m-k-1)\\& &-(m-1)(m+1)(m-3k)+3(m-1)k(m-3k)+2k(m-2)(2k-1)]\\&=&{1\over 2k(k-1)(2k-1)}\left(\!\!\!\begin{array}{c}m-k-2\\ 2k-3 \end{array}\!\!\!\right)[(m+1)(2m(m-k-1)-(m-1)(m-3k))\\& &+k(-6m(m-k-1)+3(m-1)(m-3k)+2(2k-1)(m-2)]
\\&=&{1\over 2k(k-1)(2k-1)}\left(\!\!\!\begin{array}{c}m-k-2\\ 2k-3 \end{array}\!\!\!\right)[(m+1)(m^2+km-m-3k)\\& &+k(-3m^2+km+m+k+4)]\\&=&{1\over 2k(k-1)(2k-1)}\left(\!\!\!\begin{array}{c}m-k-2\\ 2k-3 \end{array}\!\!\!\right)[(m+1)(m^2+km-m-3k)\\& &+k(m+1)(k+4-3m)]\\&=&{1\over 2k(k-1)(2k-1)}\left(\!\!\!\begin{array}{c}m-k-2\\ 2k-3 \end{array}\!\!\!\right)(m+1)(m^2+km-m-3k+k^2+4k-3km)\\&=&{1\over 2k(k-1)(2k-1)}\left(\!\!\!\begin{array}{c}m-k-2\\ 2k-3 \end{array}\!\!\!\right)(m+1)(m^2+k^2-2km-m+k)\\&=&{(m+1)(m-k)(m-k-1)\over 2k(k-1)(2k-1)}\left(\!\!\!\begin{array}{c}m-k-2\\ 2k-3 \end{array}\!\!\!\right)\\&=&{m+1\over k}\left(\!\!\!\begin{array}{c}m-k\\ 2k-1 \end{array}\!\!\!\right).\hspace{10.5cm}(2.19)\end{eqnarray*}

When $n=1$, (2.18) holds by (2.6). Let $m\in\Bbb{Z}_+$. Suppose that (2.18) holds for $n\leq m$. By (2.9),
\begin{eqnarray*}h_{m+1}&=&2Eh_m-E^2h_{m-1}+\Dlt h_{m-2}\\&=& 2E^{m+1}+[2m(m-2)-(m-1)(m-3)+2]\Dlt E^{m-2}+2\sum_{i=2}^{[|m/3|]}{m\over i}\left(\!\!\!\begin{array}{c}m-i-1\\2i-1\end{array}\!\!\!\right)\hspace{3cm}\end{eqnarray*}
\begin{eqnarray*}& &\times
\Dlt^iE^{m+1-3i}-\sum_{i=2}^{[|(m-1)/3|]}{m-1\over i}\left(\!\!\!\begin{array}{c}m-i-2\\ 2i-1 \end{array}\!\!\!\right)\Dlt^iE^{m+1-3i}\\& &+\Dlt \sum_{i=1}^{[|(m-2)/3|]}{m-2\over i}\left(\!\!\!\begin{array}{c}m-i-3\\ 2i-1\end{array}\!\!\!\right)\Dlt^iE^{m-2-3i}\\&=&
2E^{m+1}+(m+1)(m-1)\Dlt E^{m-2}+2\sum_{i=2}^{[|m/3|]}{m\over i}\left(\!\!\!\begin{array}{c}m-i-1\\ 2i-1\end{array}\!\!\!\right)\Dlt^iE^{m+1-3i}\\&&-\sum_{i=2}^{[|(m-1)/3|]}{m-1\over i}\left(\!\!\!\begin{array}{c}m-i-2\\ 2i-1\end{array}\!\!\!\right)\Dlt^iE^{m+1-3i}\\& &+\Dlt \sum_{i=2}^{[|(m-2)/3|]+1}{m-2\over i-1}\left(\!\!\!\begin{array}{c}m-i-2\\ 2i-3\end{array}\!\!\!\right)\Dlt^iE^{m+1-3i}.\hspace{5.9cm}(2.20)\end{eqnarray*}

{\it Case 1}. $m=3l$ with $l\in\Bbb{Z}_+$.
\begin{eqnarray*}h_{m+1}&=&
2E^{m+1}+(m+1)(m-1)\Dlt E^{m-2}+2\sum_{i=2}^l{m\over i}\left(\!\!\!\begin{array}{c}m-i-1\\ 2i-1\end{array}\!\!\!\right)\Dlt^iE^{m+1-3i}\\&&-\sum_{i=2}^{l-1}{m-1\over i}\left(\!\!\!\begin{array}{c}m-i-2\\ 2i-1\end{array}\!\!\!\right)\Dlt^iE^{m+1-3i}+\Dlt \sum_{i=2}^l{m-2\over i-1}\left(\!\!\!\begin{array}{c}m-i-2\\ 2i-3\end{array}\!\!\!\right)\Dlt^iE^{m+1-3i}\\&=&
2E^{m+1}+(m+1)(m-1)\Dlt E^{m-2}+\sum_{i=2}^{l-1}[{2m\over i}\left(\!\!\!\begin{array}{c}m-i-1\\ 2i-1\end{array}\!\!\!\right)-{m-1\over i}\left(\!\!\!\begin{array}{c}m-i-2\\ 2i-3\end{array}\!\!\!\right)\\& &+{m-2\over i-1}\left(\!\!\!\begin{array}{c}m-i-2\\ 2i-3\end{array}\!\!\!\right)]\Dlt^iE^{m+1-3i}+\left[{2m\over l}\left(\!\!\!\begin{array}{c}2l-1\\ 2l-1\end{array}\!\!\!\right)+{m-2\over l-2}\left(\!\!\!\begin{array}{c}2l-2\\ 2l-3\end{array}\!\!\!\right)\right]\Dlt^lE\\&=&2E^{m+1}+(m+1)(m-1)\Dlt E^{m-2}+\sum_{i=2}^{l-1}{m+1\over i}\left(\!\!\!\begin{array}{c}m-i\\ 2i-1\end{array}\!\!\!\right)\Dlt^iE^{m+1-3i}\\& &+2(m+1)\Dlt^lE\\&=&2E^{m+1}+\sum_{i=1}^{l-1}{m+1\over i}\left(\!\!\!\begin{array}{c}m-i\\ 2i-1\end{array}\!\!\!\right)\Dlt^iE^{m+1-3i}+{m+1\over l}\left(\!\!\!\begin{array}{c}m-l\\ 2l-1\end{array}\!\!\!\right)\Dlt^lE\\&=&2E^{m+1}+\sum_{i=1}^l{m+1\over i}\left(\!\!\!\begin{array}{c}m-i\\ 2i-1\end{array}\!\!\!\right)\Dlt^iE^{m+1-3i}\\&=&2E^{m+1}+\sum_{i=1}^{[|(m+1)/3|]}\left(\!\!\!\begin{array}{c}m-i\\ 2i-1\end{array}\!\!\!\right)\Dlt^iE^{m+1-3i}.\hspace{6.2cm}(2.21)\end{eqnarray*}

{\it Case 2}. $m=3l+1$ with $l\in\Bbb{N}$.
\begin{eqnarray*}h_{m+1}&=&2E^{m+1}+(m+1)(m-1)\Dlt E^{m-2}+\sum_{i=2}^l[{2m\over i}\left(\!\!\!\begin{array}{c}m-i-1\\ 2i-1\end{array}\!\!\!\right)\\& &-{m-1\over i}\left(\!\!\!\begin{array}{c}m-i-2\\ 2i-1\end{array}\!\!\!\right)+{m-2\over i-1}\left(\!\!\!\begin{array}{c}m-i-2\\ 2i-3\end{array}\!\!\!\right)]\Dlt^iE^{m+1-3i}\hspace{4cm}\end{eqnarray*}
\begin{eqnarray*}&=&2E^{m+1}+\sum_{i=1}^{[|(m+1)/3|]}{m+1\over i}\left(\!\!\!\begin{array}{c}m-i\\ 2i-1\end{array}\!\!\!\right)\Dlt^iE^{m+1-3i}.\hspace{5.5cm}(2.22)\end{eqnarray*}

{\it Case 3}. $m=3l+2$ with $l\in\Bbb{N}$.
\begin{eqnarray*}h_{m+1}&=&2E^{m+1}+(m+1)(m-1)\Dlt E^{m-2}+\sum_{i=2}^l[{2m\over i}\left(\!\!\!\begin{array}{c}m-i-1\\ 2i-1\end{array}\!\!\!\right)-{m-1\over i}\left(\!\!\!\begin{array}{c}m-i-2\\ 2i-1\end{array}\!\!\!\right)\\& &+{m-2\over (i-1)}\left(\!\!\!\begin{array}{c}m-i-2\\ 2i-3\end{array}\!\!\!\right)]\Dlt^iE^{m+1-3i}+{m-2\over l}\left(\!\!\!\begin{array}{c}2l-1\\ 2l-1\end{array}\!\!\!\right)\Dlt^{l+1}\\&=&2E^{m+1}+\sum_{i=1}^l{(m+1)\over i}\left(\!\!\!\begin{array}{c}m-i\\ 2i-1\end{array}\!\!\!\right)\Dlt^iE^{m+1-3i}+3\Dlt^{l+1}\\&=&2E^{m+1}+\sum_{i=1}^l{(m+1)\over i}\left(\!\!\!\begin{array}{c}m-i\\ 2i-1\end{array}\!\!\!\right)\Dlt^iE^{m+1-3i}+{m+1\over l+1}\left(\!\!\!\begin{array}{c}m-(l+1)\\ 2(l+1)-1\end{array}\!\!\!\right)\Dlt^{l+1}\\&=&2E^{m+1}+\sum_{i=1}^{l+1}{(m+1)\over i}\left(\!\!\!\begin{array}{c}m-i\\ 2i-1\end{array}\!\!\!\right)\Dlt^iE^{m+1-3i}\\&=&2E^{m+1}+\sum_{i=1}^{[|(m+1)/3|]}{(m+1)\over i}\left(\!\!\!\begin{array}{c}m-i\\ 2i-1\end{array}\!\!\!\right)\Dlt^iE^{m+1-3i}.\hspace{4.6cm}(2.23)\end{eqnarray*}
\pse

Thus (2.18) holds for $n=m+1$. Therefore, (2.18) holds for any $n\in\Bbb{Z}_+$ by induction on $n.\qquad\Box$
\psp

Next we consider $\rho_n$ (cf. (2.3)). Note that
$$\rho_0={(a+b)^3-a^3-b^3\over ab(a+b)}={3a^2b+3ab^2\over ab(a+b)}=3,\eqno(2.24)$$
\begin{eqnarray*}\hspace{2cm}\rho_1&=&{(a+b)^5-a^5-b^5\over ab(a+b)}\\&=&(ab)^{-1}[(a+b)^4-a^4-b^4+ab(a^2+b^2)-a^2b^2]\\&=&
(ab)^{-1}[4ab(a^2+b^2)+6a^2b^2+ab(a^2+b^2)-a^2b^2]\\&=&5(a^2+b^2+ab)\\&=&5E,\hspace{10.8cm}(2.25)\end{eqnarray*}
\begin{eqnarray*}\rho_2&=&{(a+b)^7-a^7-b^7\over ab(a+b)}\\&=&(ab)^{-1}[(a+b)^6-a^6-b^6+ab(a^4+b^4)-a^2b^2(a+b)+a^3b^3]\\&=&
(ab)^{-1}[6ab(a^4+b^4)+15a^2b^2(a+b)+20a^3b^3+ab(a^4+b^4)-a^2b^2(a+b)+a^3b^3]
\\&=&7(a^4+b^4+2ab(a+b)+3a^2b^2)\\&=&7E^2\hspace{12.8cm}(2.26)\end{eqnarray*}
by (2.4) and (2.7).
\psp

{\bf Lemma 2.3}. {\it For} $3\leq n\in\Bbb{Z}$, {\it we have}:
$$\rho_n=2E\rho_{n-1}-E^2\rho_{n-2}+\Dlt \rho_{n-3}.\eqno(2.27)$$
\pse

{\it Proof}. Note that
\begin{eqnarray*}2E\rho_{n-1}&=&(ab(a+b))^{-1}(a^2+b^2+(a+b)^2)[(a+b)^{2n+1}-a^{2n+1}-b^{2n+1}]\\&=&(ab(a+b))^{-1}[(a^2+b^2)(a+b)^{2n+1}
-a^{2n+3}-b^{2n+3}\\& &-a^2b^2(a^{2n-1}+b^{2n-1})+(a+b)^{2n+3}-(a+b)^2(a^{2n+1}+b^{2n+1})]\\&=&\rho_n+(ab(a+b))^{-1}[(a^2+b^2)(a+b)^{2n+1}-a^2b^2(a^{2n-1}+b^{2n-1})\\& &-(a+b)^2(a^{2n+1}+b^{2n+1})],\hspace{7.9cm}(2.28)\end{eqnarray*}
\begin{eqnarray*}& &E^2\rho_{n-2}\\&=&(ab(a+b))^{-1}((a+b)^2(a^2+b^2)+a^2b^2)((a+b)^{2n-1}-a^{2n-1}-b^{2n-1})\\&=&(ab(a+b))^{-1}[(a+b)^{2n+1}(a^2+b^2)-(a+b)^2(a^{2n+1}+b^{2n+1}+a^2b^2(a^{2n-3}+b^{2n-3}))\\& &+a^2b^2(a+b)^{2n-1}-a^2b^2(a^{2n-1}+b^{2n-1})]\\&=&(ab(a+b))^{-1}[(a+b)^{2n+1}(a^2+b^2)-(a+b)^2(a^{2n+1}+b^{2n+1})-\Dlt(a^{2n-3}+b^{2n-3})\\& &+\Dlt(a+b)^{2n-3}-a^2b^2(a^{2n-1}+b^{2n-1})]\
\\&=&(ab(a+b))^{-1}[(a+b)^{2n+1}(a^2+b^2)-(a+b)^2(a^{2n+1}+b^{2n+1})\\& &-a^2b^2(a^{2n-1}+b^{2n-1})]+\Dlt\rho_{n-3}
\hspace{8.4cm}(2.29)\end{eqnarray*}
by (2.11). Thus we have
$$2E\rho_{n-1}-E^2\rho_{n-2}=\rho_n-\Dlt \rho_{n-3},\eqno(2.30)$$
which is equivalent to (2.27). $\qquad\Box$
\psp

By Lemma 2.3, we calculate
$$\rho_3=9E^3+3\Dlt,\;\;\rho_4=11E^4+11\Dlt E,\;\;\rho_5=13E^5+26\Dlt E^2,\eqno(2.31)$$
$$\rho_6=15E^6+50\Dlt E^3+3\Dlt^2,\;\;\rho_7=17E^7+85\Dlt E^4+17\Dlt^2 E,\eqno(2.32)$$
$$\rho_8=19E^8+133\Dlt E^5+57\Dlt^2E^2,\;\;\rho_9=21E^9+196\Dlt E^6+147\Dlt^2 E^3+3\Dlt^3,\eqno(2.33)$$
$$\rho_{10}=23E^{10}+276\Dlt E^7+322\Dlt^2E^4+23\Dlt^3 E.\eqno(2.34)$$
Analyzing the coefficients in (2.31)-(2.34), we speculated the following theorem.
\psp

{\bf Theorem 2.4}. {\it For} $n\in\Bbb{Z}_+$, {\it we have}
$$\rho_n=\sum_{i=0}^{[|n/3|]}{2n+3\over 2i+1}\left(\!\!\!\begin{array}{c}n-i\\ 2i\end{array}\!\!\!\right)\Dlt^iE^{n-3i},\eqno(2.35)$$
{\it which is equivalent to (1.16) by (2.4) and (2.5)}.
\pse

{\it Proof}. We shall prove (2.35) by (2.27) and induction on $n$. First for $m,k\in\Bbb{Z}_+$ such that $k\leq m-1$, we have
\begin{eqnarray*}& &{2(2m+3)\over 2k+1}\left(\!\!\!\begin{array}{c}m-k\\ 2k\end{array}\!\!\!\right)-{2m+1\over 2k+1}\left(\!\!\!\begin{array}{c}m-k-1\\ 2k\end{array}\!\!\!\right)+{2m-1\over 2k-1}\left(\!\!\!\begin{array}{c}m-k-1\\ 2k-2\end{array}\!\!\!\right)\\&=&{1\over 2k(2k+1)(2k-1)}\left(\!\!\!\begin{array}{c}m-k-1\\ 2k-2\end{array}\!\!\!\right)
[2(2m+3)(m-k)(m-3k+1)\\& &-(2m+1)(m-3k+1)(m-3k)+2k(2k+1)(2m-1)]\\&=&{1\over 2k(2k+1)(2k-1)}\left(\!\!\!\begin{array}{c}m-k-1\\ 2k-2\end{array}\!\!\!\right)[2(2m+3)(m-k)(m-3k+1)\\& &-(2m+1)(m-3k+1)(m-k)+2k(2m+1)(m-3k+1)+2k(2k+1)(2m-1)]\\&=&{1\over 2k(2k+1)(2k-1)}\left(\!\!\!\begin{array}{c}m-k-1\\ 2k-2\end{array}\!\!\!\right)[(2(2m+3)-(2m+1))(m-3k+1)(m-k)\\& &+2k((2m+1)(m-3k+1)+(2k+1)(2m-1))]\\&=&{1\over 2k(2k+1)(2k-1)}\left(\!\!\!\begin{array}{c}m-k-1\\ 2k-2\end{array}\!\!\!\right)[(2m+5)(m-3k+1)(m-k)\\& &+2k(2m^2-6km+2m+m-3k+1+4km-2k+2m-1)]\\&=&{1\over 2k(2k+1)(2k-1)}\left(\!\!\!\begin{array}{c}m-k-1\\ 2k-2\end{array}\!\!\!\right)[(2m+5)(m-3k+1)(m-k)\\& &+2k(2m^2-2km+5m-5k)]\\&=&{1\over 2k(2k+1)(2k-1)}\left(\!\!\!\begin{array}{c}m-k-1\\ 2k-2\end{array}\!\!\!\right)[(2m+5)(m-3k+1)(m-k)\\& & +2k(2m+5)(m-k)]\\&=&{(2m+5)(m-k+1)(m-k)\over 2k(2k+1)(2k-1)}\left(\!\!\!\begin{array}{c}m-k-1\\ 2k-2\end{array}\!\!\!\right)\\&=&{2m+5\over 2k+1}\left(\!\!\!\begin{array}{c}m+1-k\\ 2k\end{array}\!\!\!\right).\hspace{10cm}(2.36)\end{eqnarray*}
Note that  when $k=0$, $(^{m-k-1}_{\;\;2k-2})=0$ and
$${2(2m+3)\over 2k+1}\left(\!\!\!\begin{array}{c}m-k\\ 2k\end{array}\!\!\!\right)-{2m+1\over 2k+1}\left(\!\!\!\begin{array}{c}m-k-1\\ 2k\end{array}\!\!\!\right)=2(2m+3)-(2m+1)=2m+5.\eqno(2.37)$$

When $n=1$, (2.35) holds by (2.25). Let $m\in\Bbb{Z}_+$. Suppose that (2.35) holds for $n\leq m$. By (2.35), we have the following cases.
\pse

{\it Case 1}. $m=3l$ with $l\in\Bbb{Z}_+$.
\begin{eqnarray*}\rho_{m+1}&=&2\sum_{i=0}^l{2m+3\over 2i+1}\left(\!\!\!\begin{array}{c}m-i\\ 2i\end{array}\!\!\!\right)\Dlt^iE^{m+1-3i}-\sum_{i=0}^{l-1}{2m+1\over 2i+1}\left(\!\!\!\begin{array}{c}m-i-1\\ 2i\end{array}\!\!\!\right)\Dlt^iE^{m+1-3i}\\& &+\sum_{i=0}^{l-1}{2m-1\over 2i+1}\left(\!\!\!\begin{array}{c}m-i-2\\ 2i\end{array}\!\!\!\right)\Dlt^{i+1}E^{m-3i-2}\\&=&2\sum_{i=0}^l{2m+3\over 2i+1}\left(\!\!\!\begin{array}{c}m-i\\ 2i\end{array}\!\!\!\right)\Dlt^iE^{m+1-3i}-\sum_{i=0}^{l-1}{2m+1\over 2i+1}\left(\!\!\!\begin{array}{c}m-i-1\\ 2i\end{array}\!\!\!\right)\Dlt^iE^{m+1-3i}\\& &+\sum_{i=0}^l{2m-1\over 2i-1}\left(\!\!\!\begin{array}{c}m-i-1\\ 2i-2\end{array}\!\!\!\right)\Dlt^iE^{m+1-3i}\\&=&
\sum_{i=0}^{l-1}\left[{2(2m+3)\over 2i+1}\left(\!\!\!\begin{array}{c}m-i\\ 2i\end{array}\!\!\!\right)-{2m+1\over 2i+1}\left(\!\!\!\begin{array}{c}m-i-1\\ 2i\end{array}\!\!\!\right)+{2m-1\over 2i-1}\left(\!\!\!\begin{array}{c}m-i-2\\ 2i-2\end{array}\!\!\!\right)\right]\\& &\times\Dlt^iE^{m+1-3i}+\left[{2(6l+3)\over 2l+1}\left(\!\!\!\begin{array}{c}2l\\ 2l\end{array}\!\!\!\right)+{6l-1\over 2l-1}\left(\!\!\!\begin{array}{c}3l-l-1\\ 2l-2\end{array}\!\!\!\right)\right]\Dlt^lE^{m+1-3l}\\&=&\sum_{i=0}^{l-1}{2m+5\over 2i+1}\left(\!\!\!\begin{array}{c}m+1-i\\ 2i\end{array}\!\!\!\right)\Dlt^iE^{m+1-3i}+(6l+5)\Dlt^lE^{m+1-3l}\\&=&\sum_{i=0}^{l-1}{2m+5\over 2i+1}\left(\!\!\!\begin{array}{c}m+1-i\\ 2i\end{array}\!\!\!\right)\Dlt^iE^{m+1-3i}+{2m+5\over 2l+1}\left(\!\!\!\begin{array}{c}m+1-l\\ 2l\end{array}\!\!\!\right)\Dlt^lE^{m+1-3l}\\&=&\sum_{i=0}^l{2m+5\over 2i+1}\left(\!\!\!\begin{array}{c}m+1-i\\ 2i\end{array}\!\!\!\right)\Dlt^iE^{m+1-3i}\\&=&\sum_{i=0}^{[|(m+1)/3|]}{2m+5\over 2i+1}\left(\!\!\!\begin{array}{c}m+1-i\\ 2i\end{array}\!\!\!\right)\Dlt^iE^{m+1-3i}\hspace{5.7cm}(2.38)\end{eqnarray*}

{\it Case 2}. $m=3l+1$ with $l\in\Bbb{N}$.
\begin{eqnarray*}\rho_{m+1}&=&2\sum_{i=0}^l{2m+3\over 2i+1}\left(\!\!\!\begin{array}{c}m-i\\ 2i\end{array}\!\!\!\right)\Dlt^iE^{m+1-3i}-\sum_{i=0}^l{2m+1\over 2i+1}\left(\!\!\!\begin{array}{c}m-i-1\\ 2i\end{array}\!\!\!\right)\Dlt^iE^{m+1-3i}\\& &+\sum_{i=0}^{l-1}{2m-1\over 2i+1}\left(\!\!\!\begin{array}{c}m-i-2\\ 2i\end{array}\!\!\!\right)\Dlt^{i+1}E^{m-3i-2}\\&=&\sum_{i=0}^l{2m+5\over 2i+1}\left(\!\!\!\begin{array}{c}m+1-i\\ 2i\end{array}\!\!\!\right)\Dlt^iE^{m+1-3i}\hspace{9cm}\end{eqnarray*}
\begin{eqnarray*}&=&\sum_{i=0}^{[|(m+1)/3|]}{2m+5\over 2i+1}\left(\!\!\!\begin{array}{c}m+1-i\\ 2i\end{array}\!\!\!\right)\Dlt^iE^{m+1-3i}
\hspace{6.4cm}(2.39)\end{eqnarray*}

{\it Case 3}. $m=3l+2$ with $l\in\Bbb{N}$.
\begin{eqnarray*}\rho_{m+1}&=&2\sum_{i=0}^l{2m+3\over 2i+1}\left(\!\!\!\begin{array}{c}m-i\\ 2i\end{array}\!\!\!\right)\Dlt^iE^{m+1-3i}-\sum_{i=0}^l{2m+1\over 2i+1}\left(\!\!\!\begin{array}{c}m-i-1\\ 2i\end{array}\!\!\!\right)\Dlt^iE^{m+1-3i}\hspace{1cm}\end{eqnarray*}
\begin{eqnarray*}& &+\sum_{i=0}^l{2m-1\over 2i+1}\left(\!\!\!\begin{array}{c}m-i-2\\ 2i\end{array}\!\!\!\right)\Dlt^{i+1}E^{m-3i-2}\\&=&2\sum_{i=0}^l{2m+3\over 2i+1}\left(\!\!\!\begin{array}{c}m-i\\ 2i\end{array}\!\!\!\right)\Dlt^iE^{m+1-3i}-\sum_{i=0}^l{2m+1\over 2i+1}\left(\!\!\!\begin{array}{c}m-i-1\\ 2i\end{array}\!\!\!\right)\Dlt^iE^{m+1-3i}\\& &+\sum_{i=0}^{l+1}{2m-1\over 2i-1}\left(\!\!\!\begin{array}{c}m-i-1\\ 2i-2\end{array}\!\!\!\right)\Dlt^iE^{m+1-3i}\\&=&\sum_{i=0}^l{2m+5\over 2i+1}\left(\!\!\!\begin{array}{c}m+1-i\\ 2i\end{array}\!\!\!\right)\Dlt^iE^{m+1-3i}+{6l+3\over 2l+1}\left(\!\!\!\begin{array}{c}2l\\ 2l\end{array}\!\!\!\right)\Dlt^{l+1}E^{m-3l-2}\\&=&\sum_{i=0}^l{2m+5\over 2i+1}\left(\!\!\!\begin{array}{c}m+1-i\\ 2i\end{array}\!\!\!\right)\Dlt^iE^{m+1-3i}+3\Dlt^{l+1}E^{m-3l-2}\\&=&\sum_{i=0}^l{2m+5\over 2i+1}\left(\!\!\!\begin{array}{c}m+1-i\\ 2i\end{array}\!\!\!\right)\Dlt^iE^{m+1-3i}\\& &+{2m+5\over 2(l+1)+1}\left(\!\!\!\begin{array}{c}m+1-(l+1)\\ 2(l+1)\end{array}\!\!\!\right)
\Dlt^{l+1}E^{m+1-3(l+1)}\\&=&\sum_{i=0}^{l+1}{2m+5\over 2i+1}\left(\!\!\!\begin{array}{c}m+1-i\\ 2i\end{array}\!\!\!\right)\Dlt^iE^{m+1-3i}
\\&=&\sum_{i=0}^{[|(m+1)/3|]}{2m+5\over 2i+1}\left(\!\!\!\begin{array}{c}m+1-i\\ 2i\end{array}\!\!\!\right)\Dlt^iE^{m+1-3i}.\hspace{6.4cm}(2.40)\end{eqnarray*}
\pse

Thus (2.35) holds for $n=m+1$. Therefore, (2.34) holds for any $n\in\Bbb{Z}_+$ by induction on $n.\qquad\Box$

\section{Theta Series and Weighted Symmetric Polynomials}

In this section, we shall determine the theta series of  certain infinite families of positive definite even unimodular lattices containing a sublattice of the same rank and isomorphic to the direct sum of finite copies of the lattices $R_{D_{2n}}$ with various $n$. The theta series of these lattices are weighted symmetric polynomials of the functions of $\{\Dlt_{24}(z),h_n(z),\rho_n(z)\mid n\in\Bbb{Z}_+\}$. By (1.15) and (1.16), we essentially determine  the theta series of these lattices as polynomials of the well-known Essenstein series $E_4(z)$ and Ramanujan series $\Dlt_{24}(z)$.

Let $k$ be a positive integer. Suppose that $\{n_1,...,n_k\}\subset \Bbb{Z}_+$ is a subset. Set
$$\bar{n}_0=0,\;\;\bar{n}_i=\sum_{j=1}^in_i,\;\;n=\bar{n}_k.\eqno(3.1)$$
 Recall the Euclidean space $\Bbb{R}^n$ with inner product (1.8). All the lattices in this section have the symmetric bilinear form inherited from this inner product.
Set
$$Q=\{\vec{\al}=(\al_1,...,\al_n)\in\Bbb{Z}^n\mid \sum_{j=1}^{n_i}\al_{\bar{n}_{i-1}+j}\in 2\Bbb{Z}\;\for\;i\in\ol{1,k}\}\subset \Bbb{R}^n.\eqno(3.2)$$
Then $Q$ is a lattice that isomorphic to $R_{D_{n_1}}\oplus\cdots\oplus R_{D_{n_k}}$. For $j\in\Bbb{N}$, we denote 
$${\bf 0}_j=(0,...,0),\;\;{\bf 1}_j=(1,...,1)\in \Bbb{R}^j.\eqno(3.3)$$ 
Set
$$x^i_1={1\over 2}({\bf 0}_{\bar{n}_{i-1}},{\bf 1}_{n_i},{\bf 0}_{n-\bar{n}_i}),\;\;x^i_2=({\bf 0}_{\bar{n}_{i-1}},1,{\bf 0}_{n-\bar{n}_{i-1}-1}),\;\;x^i_3=x^i_1+x^i_2\in\Bbb{R}^n\eqno(3.4)$$
for $i=1,...,k$. When $k=1$, we have
$$\Theta_Q(z)=\Theta_{R_{D_{n_1}}}={1\over 2}(\vt_3(z)^{n_1}+\vt_4(z)^{n_1}),\qquad\Theta_{x^1_1+Q}(z)=\Theta_{x^1_3+Q}(z)={1\over 2}\vt_2(z)^{n_1},\eqno(3.5)$$
$$\Theta_{x^1_2+Q}(z)={1\over 2}(\vt_3(z)^{n_1}-\vt_4(z)^{n_1})\eqno(3.6)$$
(cf. (1.12) and Chapter 4 of [CS3]). Moreover,
 $$ 2x^1_1,2x^1_2\in Q\eqno(3.7)$$
when $n_1$ is even. In fact, we have
$$\{\vec{\al}\in\Bbb{R}^n\mid \la\vec{\al},\vec{\be}\ra \in\Bbb{Z}\;\for\;\vec{\be}\in Q\}=\sum_{i=1}^k(\Bbb{Z}x^i_1+\Bbb{Z}x^i_2)+Q\eqno(3.8)$$
(eg., cf. [CS3] or [X3]).

Set
$$L=\sum_{i=1}^k\Bbb{Z}(x_1^i+x^1_2+\cdots +x^{i-1}_2+x^{i+1}_2+\cdots +x^k_2)+Q.\eqno(3.9)$$
When $k=1$,
$$L=\Bbb{Z}x^1_1+Q.\eqno(3.10)$$
\psp

{\bf Theorem 3.1}. {\it If} $k=2\ell+1$ {\it is an odd positive integer and} $n_i=8m_i$ {\it with} $\ell\in\Bbb{N}$ {\it and} $m_i\in\Bbb{Z}_+$, {\it then}
$L$ {\it is a positive definite even unimodular lattice with the theta series}
\begin{eqnarray*}& &2^{2\ell+1}\Theta_L(z)\\&=&\sum_{1\leq j_1\leq j_2;j_1+j_2\leq \ell}\sym\{h_{m_1+\cdots +m_{2j_1}}(z)h_{m_{2j_1+1}+\cdots+m_{2(j_1+j_2)}}(z)h_{m_{2(j_1+j_2)+1}+\cdots+m_{2\ell+1}}(z)\}\hspace{2cm}\end{eqnarray*}
\begin{eqnarray*}
& &-\sum_{0\leq j_1\leq j_2\leq \ell-j_1-j_2-1}\sym\{h_{m_1+\cdots+m_{2j_1+1}}(z)h_{m_{2j_1+2}+\cdots+m_{2(j_1+j_2+1)}}(z)\\& &\times h_{m_{2(j_1+j_2)+3}+\cdots+m_{2\ell+1}}(z)\}
+\sum_{j=1}^{\ell}(3-4^{\ell-j})\sym\{h_{m_1+\cdots +m_{2j}}(z) h_{m_{2j+1}+\cdots +m_{2\ell+1}}(z)\}\\&&+\left(4-3\cdot 4^{\ell}+\frac{2}{3}\sum_{j_1+j_2\leq \ell-1}\left(\!\!\!\begin{array}{c}2\ell+1\\ 2j_1+1,2j_2+1\end{array}\!\!\!\right)\right)h_{m_1+\cdots+ m_{2\ell+1}}(z)\hspace{3.1cm}(3.11)\end{eqnarray*}

{\it If} $k=2\ell$ {\it is an even positive integer and} $n_i=8m_i+4$ {\it with} $\ell\in\Bbb{Z}_+$ {\it and} $m_i\in\Bbb{N}$, {\it then}
$L$ {\it is a positive definite even unimodular lattice} {\it with the theta series}
\begin{eqnarray*}& &2^{2\ell}\Theta_L(z)\\&=&\sum_{1\leq j_1\leq j_2\leq\ell-j_1-j_2}\sym\{h_{m_1+\cdots+m_{2j_1}+j_1}(z)h_{m_{2j_1+1}+\cdots+m_{2(j_1+j_2)}+j_2}(z)\\& &\times h_{m_{2(j_1+j_2)+1}+\cdots+m_{2\ell}+\ell-j_1-j_2}(z)\}
-2^8\Dlt_{24}(z)\sum_{1\leq j_1,j_2;j_1+j_2\leq\ell-2}\sym\{\rho_{m_1+\cdots+m_{2j_1+1}+j_1-1}(z)\\& &\times\rho_{m_{2j_1+2}+\cdots+m_{2(j_1+j_2+1)}+j_2-1}(z)h_{m_{2(j_1+j_2)+3}+\cdots+m_{2\ell}+\ell-j_1-j_2-1}(z)\}\\& &+3\sum_{j=1}^{[|\ell/2|]}\sym\{h_{m_1+\cdots+m_{2j}+j}(z)
 h_{m_{2j+1}+\cdots+m_{2\ell}+\ell-j}(z)\}+2^8\Dlt_{24}(z)\\& &\times\sum_{j=0}^{[|(\ell-1)/2|]}(2^{2j}+2^{2(\ell-j-1)}-3)\sym\{\rho_{m_1+\cdots+m_{2j+1}+j-1}(z)\rho_{m_{2j+2}+\cdots+m_{2\ell}+\ell-j-2}(z)\}\\& &+\left(4-\frac{2}{3}\sum_{j_1+j_2\leq\ell}\left(\!\!\!\begin{array}{c}2\ell\\ 2j_1,2j_2\end{array}\!\!\!\right)\right)h_{m_1+\cdots+m_{2\ell}+\ell}(z)\hspace{5.9cm}(3.12)\end{eqnarray*}
\pse

{\it Proof}. First we consider the case that $k=2\ell+1$ {\it is an odd positive number and} $n_i=8m_i$ {\it with} $\ell\in\Bbb{N}$ {\it and} $m_i\in\Bbb{Z}_+$. The lattice $L$ is unimodular by (1.11) and (3.7)-(3.9). Before we prove (3.11), we need some combinatorial facts as follows. For $j\in\Bbb{Z}_+$, we have
$$\sum_{i=0}^{[|j/2|]}\left(\!\!\!\begin{array}{c}j\\ 2i\end{array}\!\!\!\right)\pm \sum_{i=0}^{[|(j-1)/2|]}\left(\!\!\!\begin{array}{c}j\\ 2i+1\end{array}\!\!\!\right)=(1\pm 1)^j,\eqno(3.13)$$
which implies 
$$\sum_{i=0}^{[|j/2|]}\left(\!\!\!\begin{array}{c}j\\ 2i\end{array}\!\!\!\right)=\sum_{i=0}^{[|(j-1)/2|]}\left(\!\!\!\begin{array}{c}j\\ 2i+1\end{array}\!\!\!\right)=2^{j-1}.\eqno(3.14)$$
Thus
$$\sum_{i=1}^{[|j/2|]}\left(\!\!\!\begin{array}{c}j\\ 2i\end{array}\!\!\!\right)=2^{j-1}-1.\eqno(3.15)$$
Moreover,
$$\sum_{i=0}^j\left(\!\!\!\begin{array}{c}j\\ 2i\end{array}\!\!\!\right)2^{j-2i}=\frac{(1+2)^j+(-1+2)^j}{2}=\frac{3^j+1}{2}.\eqno(3.16)$$
In particular,
$$\sum_{i=1}^{\ell}\left(\!\!\!\begin{array}{c}2\ell+1\\ 2i\end{array}\!\!\!\right)4^{\ell-i}=\frac{3^{2\ell+1}+1}{4}-4^{\ell}.\eqno(3.17)$$
Furthermore,
\begin{eqnarray*}\sum_{j_1+j_2\leq \ell}\left(\!\!\!\begin{array}{c}2\ell+1\\ 2j_1,2j_2\end{array}\!\!\!\right)+\sum_{j_1+j_2\leq \ell-1}\left(\!\!\!\begin{array}{c}2\ell+1\\ 2j_1+1,2j_2+1\end{array}\!\!\!\right)&=&\frac{(1+1+1)^{2\ell+1}+(-1-1+1)^{2\ell+1}}{2}\\&=&\frac{3^{2\ell+1}-1}{2}.\hspace{3.8cm}(3.18)\end{eqnarray*}

The following countings are also needed. For $\ell>j\in\Bbb{Z}_+$, we have
\begin{eqnarray*}& &\mbox{the number of minomials in the expression}
\sum_{i=0}^{[|(j-1)/2|]}\sym\{h_{m_1+\cdots+m_{2i+1}}h_{m_{2i+2}+\cdot+m_{2j}}\}\\&=&{1\over 2}\sum_{i=0}^{j-1}\left(\!\!\!\begin{array}{c}2j\\ 2i+1\end{array}\!\!\!\right),\hspace{11cm}(3.19)\end{eqnarray*}
\begin{eqnarray*}& &\mbox{the number of minomials in the expression}\sum_{i=1}^{[|j/2|]}\sym\{h_{m_1+\cdots+m_{2i}}h_{m_{2i+1}+\cdot+m_{2j}}\}\\&=&{1\over 2}\sum_{i=1}^{j-1}\left(\!\!\!\begin{array}{c}2j\\ 2i\end{array}\!\!\!\right).\hspace{11.7cm}(3.20)\end{eqnarray*} 
Moreover,
\begin{eqnarray*}& &\mbox{the number of  minomials in the expression}\\& &
 \sum_{1\leq j_1\leq j_2;j_1+j_2\leq \ell}\sym\{h_{m_1+\cdots +m_{2j_1}}(z)h_{m_{2j_1+1}+\cdots+m_{2(j_1+j_2)}}(z)h_{m_{2(j_1+j_2)+1}+\cdots+m_{2\ell+1}}(z)\}\\&=&\frac{1}{2}\sum_{1\leq j_1,j_2;j_1+j_2\leq \ell}\left(\!\!\!\begin{array}{c}2\ell+1\\ 2j_1,2j_2\end{array}\!\!\!\right),\hspace{9.2cm}(3.21)\end{eqnarray*}
\begin{eqnarray*}& &\mbox{the number of  minomials in the expression}\;\;\sum_{0\leq j_1\leq j_2\leq \ell-j_1-j_2-1}\\& &\sym\{h_{m_1+\cdots+m_{2j_1+1}}(z)h_{m_{2j_1+2}+\cdots+m_{2(j_1+j_2+1)}}(z)h_{m_{2(j_1+j_2)+3}+\cdots+m_{2\ell+1}}(z)\}\\&=&\frac{1}{6}\sum_{j_1+j_2\leq \ell-1}\left(\!\!\!\begin{array}{c}2\ell+1\\ 2j_1+1,2j_2+1\end{array}\!\!\!\right).\hspace{8.5cm}(3.22)\end{eqnarray*}
Furthermore,
\begin{eqnarray*}\hspace{2cm}& &\sum_{j_1+j_2\leq \ell}\left(\!\!\!\begin{array}{c}2\ell+1\\ 2j_1,2j_2\end{array}\!\!\!\right)\\&=&2\sum_{i=0}^{\ell}\left(\!\!\!\begin{array}{c}2\ell+1\\ 2i\end{array}\!\!\!\right)-1+\sum_{1\leq j_1\leq j_2;j_1+j_2\leq \ell}\left(\!\!\!\begin{array}{c}2\ell+1\\ 2j_1,2j_2\end{array}\!\!\!\right)\\&=& 2^{2\ell+1}-1+\sum_{1\leq j_1\leq j_2;j_1+j_2\leq \ell}\left(\!\!\!\begin{array}{c}2\ell+1\\ 2j_1,2j_2\end{array}\!\!\!\right).\hspace{5.4cm}(3.23)\end{eqnarray*}

By (1.12), (3.5)-(3.7), (3.13)-(3.23) and calculating the theta series of the following cosets of $Q$:
$$Q,\qquad\sum_{i=1}^{2j}(x^i_1+x^1_2+\cdots +x^{i-1}_2+x^{i+1}_2+\cdots +x^k_2)+Q,\eqno(3.24)$$
and
$$\sum_{i=1}^{2j'+1}(x^i+x^1_2+\cdots +x^{i-1}_2+x^{i+1}_2+\cdots +x^k_2)+Q\eqno(3.25)$$
for $1\leq j\leq \ell$ and $0\leq j'\leq \ell$, we obtain
\begin{eqnarray*}& &2^{2\ell+1}\Theta_L(z)\\&=&\sum_{j=0}^{\ell}\sym\{\vt_2(z)^{n_1+\cdots +n_{2j}}\prod_{p=2j+1}^{2\ell+1}(\vt_3(z)^{n_p}+\vt_4(z)^{n_p})\}\\& &+\sum_{j=0}^{\ell}\sym\{\vt_2(z)^{n_1+\cdots +n_{2j+1}}\prod_{p=2j+2}^{2\ell+1}(\vt_3(z)^{n_p}-\vt_4(z)^{n_p})\}\\&=& \sum_{1\leq j_1\leq j_2;j_1+j_2\leq \ell}[\sym\{\vt_2(z)^{n_1+\cdots +n_{2j_1}}\vt_3(z)^{n_{2j_1+1}+\cdots +n_{2(j_1+j_2)}}\vt_4(z)^{n_{2(j_1+j_2)+1}+\cdots+n_{2\ell+1}}\}
\\& &+\sym\{\vt_3(z)^{n_1+\cdots +n_{2j_1}}\vt_4(z)^{n_{2j_1+1}+\cdots +n_{2(j_1+j_2)}}\vt_2(z)^{n_{2(j_1+j_2)+1}+\cdots+ n_{2\ell+1}}\}\\& &+\sym\{\vt_4(z)^{n_1+\cdots +n_{2j_1}}\vt_2(z)^{n_{2j_1+1}+\cdots +n_{2(j_1+j_2)}}\vt_3(z)^{n_{2(j_1+j_2)+1}+\cdots +n_{2\ell+1}}\}]
\\& &-\sum_{0\leq j_1\leq j_2\leq \ell-j_1-j_2-1}\sym\{\vt_2(z)^{n_1+\cdots +n_{2j_1+1}}
\vt_3(z)^{n_{2j_1+2}+\cdots +n_{2(j_1+j_2+1)}}\\& &\times\vt_4(z)^{n_{2(j_1+j_2)+3}+\cdots +n_{2\ell+1}}\}+\sum_{j=1}^{\ell}[\sym\{
\vt_2(z)^{n_1+\cdots +n_{2j}}(\vt_3(z)^{n_{2j+1}+\cdots +n_{2\ell+1}}\\& &+\vt_4(z)^{n_{2j+1}+\cdots +n_{2\ell+1}})\}+\sym\{
\vt_3(z)^{n_1+\cdots +n_{2j}}(\vt_2(z)^{n_{2j+1}+\cdots +n_{2\ell+1}}\\& &+\vt_4(z)^{n_{2j+1}+\cdots +n_{2\ell+1}})\}+\sym\{\vt_4(z)^{n_1+\cdots +n_{2j}}(\vt_3(z)^{n_{2j+1}+\cdots +n_{2\ell+1}}\\&&+\vt_2(z)^{n_{2j+1}+\cdots +n_{2\ell+1}})\}]+\vt_2(z)^{n_1+\cdots+n_{2\ell+1}}+\vt_3(z)^{n_1+\cdots+n_{2\ell+1}}+\vt_4(z)^{n_1+\cdots+n_{2\ell+1}}\\&=&\sum_{1\leq j_1\leq j_2;j_1+j_2\leq \ell}\sym\{(\vt_2(z)^{n_1+\cdots +n_{2j_1}}+\vt_3(z)^{n_1+\cdots +n_{2j_1}}+\vt_4(z)^{n_1+\cdots +n_{2j_1}})\\& &\times(\vt_2(z)^{n_{2j_1+1}+\cdots +n_{2(j_1+j_2)}}+\vt_3(z)^{n_{2j_1+1}+\cdots +n_{2(j_1+j_2)}}+\vt_4(z)^{n_{2j_1+1}+\cdots +n_{2(j_1+j_2)}})\\& &\times(\vt_2(z)^{n_{2(j_1+j_2)+1}+\cdots+n_{2\ell+1}}+\vt_3(z)^{n_{2(j_1+j_2)+1}+\cdots+n_{2\ell+1}}+\vt_4(z)^{n_{2(j_1+j_2)+1}+\cdots+n_{2\ell+1}})\}
\\& &-\sum_{0\leq j_1\leq j_2\leq \ell-j_1-j_2-1}\sym\{(\vt_2(z)^{n_1+\cdots +n_{2j_1+1}}+\vt_3(z)^{n_1+\cdots +n_{2j_1+1}}+\vt_4(z)^{n_1+\cdots +n_{2j_1+1}})\\&&\times (\vt_2(z)^{n_{2j_1+2}+\cdots +n_{2(j_1+j_2+1)}}+\vt_3(z)^{n_{2j_1+2}+\cdots +n_{2(j_1+j_2+1)}}+\vt_4(z)^{n_{2j_1+2}+\cdots +n_{2(j_1+j_2+1)}})\\& &\times(\vt_2(z)^{n_{2(j_1+j_2)+3}+\cdots+ n_{2\ell+1}}+\vt_3(z)^{n_{2(j_1+j_2)+3}+\cdots+n_{2\ell+1}}+\vt_4(z)^{n_{2(j_1+j_2)+3}+\cdots+n_{2\ell+1}})\}\\& &+\sum_{j=1}^{\ell}\left(1+\frac{1}{2}\sum_{i=0}^{j-1}\left(\!\!\!\begin{array}{c}2j\\ 2i+1\end{array}\!\!\!\right)-\frac{1}{2}\sum_{i=1}^{j-1}\left(\!\!\!\begin{array}{c}2j\\ 2i\end{array}\!\!\!\right)-\sum_{i=1}^{\ell-j}\left(\!\!\!\begin{array}{c}2(\ell-j)+1\\ 2i\end{array}\!\!\!\right)\right)\\& &\times [\sym\{
\vt_2(z)^{n_1+\cdots +n_{2j}}(\vt_3(z)^{n_{2j+1}+\cdots +n_{2\ell+1}}+\vt_4(z)^{n_{2j+1}+\cdots +n_{2\ell+1}})\}\hspace{6cm}\end{eqnarray*}
\begin{eqnarray*}& &+\sym\{
\vt_3(z)^{n_1+\cdots +n_{2j}}(\vt_2(z)^{n_{2j+1}+\cdots +n_{2\ell+1}}+\vt_4(z)^{n_{2j+1}+\cdots +n_{2\ell+1}})\}\\& &+\sym\{\vt_4(z)^{n_1+\cdots +n_{2j}}(\vt_3(z)^{n_{2j+1}+\cdots +n_{2\ell+1}}+\vt_4(z)^{n_{2j+1}+\cdots +n_{2\ell+1}})]
\\& &+\left(1-\frac{1}{2}\sum_{1\leq j_1,j_2;j_1+j_2\leq \ell}\left(\!\!\!\begin{array}{c}2\ell+1\\ 2j_1,2j_2\end{array}\!\!\!\right)+\frac{1}{6}\sum_{j_1+j_2\leq\ell-1}\left(\!\!\!\begin{array}{c}2\ell+1\\ 2j_1+1,2j_2+1\end{array}\!\!\!\right)\right)\\& &\times(\vt_2(z)^{n_1+\cdots+n_{2\ell+1}}+\vt_3(z)^{n_1+\cdots+n_{2\ell+1}}+\vt_4(z)^{n_1+\cdots+n_{2\ell+1}})\\&=&
\sum_{1\leq j_1\leq j_2;j_1+j_2\leq \ell}\sym\{h_{m_1+\cdots+m_{2j_1}}(z)h_{m_{2j_1+1}+\cdots+m_{2(j_1+j_2)}}(z)h_{m_{2(j_1+j_2)+1}+\cdots+m_{2\ell+1}}(z)\}\\& &-\sum_{0\leq j_1\leq j_2\leq \ell-j_1-j_2-1}\sym\{h_{m_1+\cdots+m_{2j_1+1}}(z)h_{m_{2j_1+2}+\cdots+m_{2(j_1+j_2+1)}}(z)\\& &\times h_{m_{2(j_1+j_2)+3}+\cdots+m_{2\ell+1}}(z)\}+\sum_{j=1}^{\ell}(3-4^{\ell-j})\sym\{(
\vt_2(z)^{n_1+\cdots +n_{2j}}+\vt_3(z)^{n_1+\cdots +n_{2j}}\\& &+\vt_4(z)^{n_1+\cdots +n_{2j}})(\vt_2(z)^{n_{2j+1}+\cdots +n_{2\ell+1}}+\vt_3(z)^{n_{2j+1}+\cdots +n_{2\ell+1}}+\vt_4(z)^{n_{2j+1}+\cdots +n_{2\ell+1}})\}\\& &+[2^{2\ell}+\frac{1}{2}-\frac{1}{2}\sum_{j_1+j_2\leq \ell}\left(\!\!\!\begin{array}{c}2\ell+1\\ 2j_1,2j_2\end{array}\!\!\!\right)+\frac{1}{6}\sum_{j_1+j_2\leq \ell-1}\left(\!\!\!\begin{array}{c}2\ell+1\\ 2j_1+1,2j_2+1\end{array}\!\!\!\right)-\sum_{j=1}^{\ell}\left(\!\!\!\begin{array}{c}2\ell+1\\ 2j\end{array}\!\!\!\right)\\& &\times(3-4^{\ell-j})](\vt_2(z)^{n_1+\cdots+n_{2\ell+1}}+\vt_3(z)^{n_1+\cdots+n_{2\ell+1}}+\vt_4(z)^{n_1+\cdots+n_{2\ell+1}})\\&=& \sum_{1\leq j_1\leq j_2;j_1+j_2\leq \ell}\sym\{h_{m_1+\cdots +m_{2j_1}}(z)h_{m_{2j_1+1}+\cdots+m_{2(j_1+j_2)}}(z)h_{m_{2(j_1+j_2)+1}+\cdots+m_{2\ell+1}}(z)\}\\& &-\sum_{0\leq j_1\leq j_2\leq \ell-j_1-j_2-1}\sym\{h_{m_1+\cdots+m_{2j_1+1}}(z)h_{m_{2j_1+2}+\cdots+m_{2(j_1+j_2+1)}}(z)\\& &\times h_{m_{2(j_1+j_2)+3}+\cdots+m_{2\ell+1}}(z)\}
+\sum_{j=1}^{\ell}(3-4^{\ell-j})\sym\{h_{m_1+\cdots +m_{2j}}(z) h_{m_{2j+1}+\cdots +m_{2\ell+1}}(z)\}\\&&+\left(4-3\cdot 4^{\ell}+\frac{2}{3}\sum_{j_1+j_2\leq \ell-1}\left(\!\!\!\begin{array}{c}2\ell+1\\ 2j_1+1,2j_2+1\end{array}\!\!\!\right)\right)h_{m_1+\cdots+ m_{2\ell+1}}(z)\hspace{2.8cm}(3.26)\end{eqnarray*}

Next we consider the case that $k=2\ell$ is an even positive number and $n_i=8m_i+4$. $L$ is unimodular again by (1.11) and (3.7)-(3.9). We need the combinatorial facts as follows. First,
\begin{eqnarray*}\hspace{2cm}& &\sum_{j=0}^{\ell-1}\left(\!\!\!\begin{array}{c}2\ell\\ 2j+1\end{array}\!\!\!\right)(2^{2j}+2^{2(\ell-j-1)})=\sum_{j=0}^{\ell-1}\left(\!\!\!\begin{array}{c}2\ell\\ 2j+1\end{array}\!\!\!\right)2^{2j+1}\\&=&\frac{(1+2)^{2\ell}-(-1+2)^{2\ell}}{2}=\frac{3^{2\ell}-1}{2}.\hspace{5.9cm}(3.27)\end{eqnarray*}
Moreover,
\begin{eqnarray*}\sum_{j_1+j_2\leq \ell}\left(\!\!\!\begin{array}{c}2\ell\\ 2j_1,2j_2\end{array}\!\!\!\right)&=&3\sum_{j=0}^{\ell}\left(\!\!\!\begin{array}{c}2\ell\\ 2j\end{array}\!\!\!\right)-3+\sum_{1\leq j_1,j_2;j_1+j_2\leq \ell-1}\left(\!\!\!\begin{array}{c}2\ell\\ 2j_1,2j_2\end{array}\!\!\!\right)\\&=& 3(2^{2\ell-1}-1)+\sum_{1\leq j_1,j_2;j_1+j_2\leq \ell-1}\left(\!\!\!\begin{array}{c}2\ell\\ 2j_1,2j_2\end{array}\!\!\!\right),\hspace{3.4cm}(3.28)\end{eqnarray*}
\begin{eqnarray*}\sum_{j_1+j_2\leq \ell-1}\left(\!\!\!\begin{array}{c}2\ell\\ 2j_1+1,2j_2+1\end{array}\!\!\!\right)&=&\sum_{j=0}^{\ell-1}\left(\!\!\!\begin{array}{c}2\ell\\ 2j+1\end{array}\!\!\!\right)+\sum_{j_1+j_2\leq \ell-2}\left(\!\!\!\begin{array}{c}2\ell\\ 2j_1+1,2j_2+1\end{array}\!\!\!\right)\\&=&2^{2\ell-1}+\sum_{j_1+j_2\leq \ell-2}\left(\!\!\!\begin{array}{c}2\ell\\ 2j_1+1,2j_2+1\end{array}\!\!\!\right)\hspace{2.7cm}(3.29)\end{eqnarray*}
by (3.14). Furthermore,
\begin{eqnarray*}\sum_{j_1+j_2\leq \ell}\left(\!\!\!\begin{array}{c}2\ell\\ 2j_1,2j_2\end{array}\!\!\!\right)+\sum_{j_1+j_2\leq \ell-1}\left(\!\!\!\begin{array}{c}2\ell\\ 2j_1+1,2j_2+1\end{array}\!\!\!\right)&=&\frac{(1+1+1)^{2\ell}+(-1-1+1)^{2\ell}}{2}\\&=&\frac{3^{2\ell}-1}{2}.\hspace{4.1cm}(3.30)\end{eqnarray*}

The following countings are also needed. Note
\begin{eqnarray*}& &\mbox{the number of minomials in the expression}\sum_{1\leq j_1\leq j_2;j_1+j_2\leq \ell-1}\\& &
 \sym\{h_{m_1+\cdots +m_{2j_1}}(z)h_{m_{2j_1+1}+\cdots+m_{2(j_1+j_2)}}(z)h_{m_{2(j_1+j_2)+1}+\cdots+m_{2\ell}}(z)\}\\&=&\frac{1}{6}\sum_{1\leq j_1,j_2;j_1+j_2\leq \ell-1}\left(\!\!\!\begin{array}{c}2\ell\\ 2j_1,2j_2\end{array}\!\!\!\right),\hspace{8.9cm}(3.31)\end{eqnarray*}
\begin{eqnarray*}& &\mbox{the number of minomials in the expression}\;\;\sum_{0\leq j_1\leq j_2;j_1+j_2\leq \ell-2}\\& &\sym\{h_{m_1+\cdots+m_{2j_1+1}}(z)h_{m_{2j_1+2}+\cdots+m_{2(j_1+j_2+1)}}(z)h_{m_{2(j_1+j_2)+3}+\cdots+m_{2\ell}}(z)\}\\&=&\frac{1}{2}\sum_{j_1+j_2\leq \ell-2}\left(\!\!\!\begin{array}{c}2\ell\\ 2j_1+1,2j_2+1\end{array}\!\!\!\right).\hspace{8.7cm}(3.32)\end{eqnarray*}

By (3.5), (3.6), (3.13)-(3.15), (3.19), (3.20), (3.27)-(3.32) and calculating the theta series of the cosets in $L/Q$ by (3.5)-(3.7), we have
\begin{eqnarray*}& &2^{2\ell+1}\Theta_L(z)\\&=&\sum_{j=0}^{\ell}\sym\{\vt_2(z)^{n_1+\cdots +n_{2j}}\prod_{p=2j+1}^{2\ell}(\vt_3(z)^{n_p}+\vt_4(z)^{n_p})\}\\& &+\sum_{j=0}^{\ell-1}\sym\{\vt_2(z)^{n_1+\cdots +n_{2j+1}}\prod_{p=2j+2}^{2\ell}(\vt_3(z)^{n_p}-\vt_4(z)^{n_p})\}\\&=&\sum_{0\leq j_1\leq j_2\leq\ell-j_1-j_2}\sym\{\vt_2(z)^{n_1+\cdots +n_{2j_1}}\vt_3(z)^{n_{2j_1+1}+\cdots+n_{2(j_1+j_2)}}\vt_4(z)^{n_{2(j_1+j_2)+1}+\cdots+n_{2\ell}}\}\\&&+\sum_{0\leq j_1\leq j_2;j_1+j_2\leq\ell-1}\sym\{\vt_2(z)^{n_1+\cdots +n_{2j_1+1}}\vt_3(z)^{n_{2j_1+2}+\cdots+n_{2(j_1+j_2+1)}}\vt_4(z)^{n_{2(j_1+j_2)+3}+\cdots+n_{2\ell}}\}\\&&-\sum_{0\leq j_1\leq j_2;j_1+j_2\leq\ell-1}\sym\{\vt_2(z)^{n_1+\cdots +n_{2j_1+1}}\vt_4(z)^{n_{2j_1+2}+\cdots+n_{2(j_1+j_2+1)}}\vt_3(z)^{n_{2(j_1+j_2)+3}+\cdots+n_{2\ell}}\}\\&&+\sum_{0\leq j_1\leq j_2;j_1+j_2\leq\ell-1}\sym\{\vt_3(z)^{n_1+\cdots +n_{2j_1+1}}\vt_4(z)^{n_{2j_1+2}+\cdots+n_{2(j_1+j_2+1)}}\vt_2(z)^{n_{2(j_1+j_2)+3}+\cdots+n_{2\ell}}\}\hspace{8cm}\end{eqnarray*}
\begin{eqnarray*}&=&
\sum_{1\leq j_1\leq j_2\leq\ell-j_1-j_2}\sym\{(\vt_2(z)^{n_1+\cdots +n_{2j_1}}+\vt_3(z)^{n_1+\cdots +n_{2j_1}}+\vt_4(z)^{n_1+\cdots +n_{2j_1}})\\& &\times (\vt_2(z)^{n_{2j_1+1}+\cdots+n_{2(j_1+j_2)}}+\vt_3(z)^{n_{2j_1+1}+\cdots+n_{2(j_1+j_2)}}+\vt_4(z)^{n_{2j_1+1}+\cdots+n_{2(j_1+j_2)}})\\& &\times(\vt_2(z)^{n_{2(j_1+j_2)+1}+\cdots+n_{2\ell}}+\vt_3(z)^{n_{2(j_1+j_2)+1}+\cdots+n_{2\ell}}+\vt_4(z)^{n_{2(j_1+j_2)+1}+\cdots+n_{2\ell}})\}
\\& &-\sum_{0 \leq j_1\leq j_2;j_1+j_2\leq\ell-2}\sym\{(\vt_3(z)^{n_1+\cdots +n_{2j_1+1}}-\vt_2(z)^{n_1+\cdots +n_{2j_1+1}}-\vt_4(z)^{n_1+\cdots +n_{2j_1+1}})\\& &\times (\vt_3(z)^{n_{2j_1+2}+\cdots+n_{2(j_1+j_2+1)}}-\vt_2(z)^{n_{2j_1+2}+\cdots+n_{2(j_1+j_2+1)}}-\vt_4(z)^{n_{2j_1+2}+\cdots+n_{2(j_1+j_2+1)}})\\& &\times(\vt_2(z)^{n_{2(j_1+j_2)+3}+\cdots+n_{2\ell}}+\vt_3(z)^{n_{2(j_1+j_2)+3}+\cdots+n_{2\ell}}+\vt_4(z)^{n_{2(j_1+j_2)+3}+\cdots+n_{2\ell}})\}
\\& &+3\sum_{j=1}^{[|\ell/2|]}\sym\{(\vt_2(z)^{n_1+\cdots+n_{2j}}+\vt_3(z)^{n_1+\cdots+n_{2j}}+\vt_4(z)^{n_1+\cdots+n_{2j}})(\vt_2(z)^{n_{2j+1}+\cdots+n_{2\ell}}\\& &+\vt_3(z)^{n_{2j+1}+\cdots+n_{2\ell}}+\vt_4(z)^{n_{2j+1}+\cdots+n_{2\ell}})\}+\sum_{j=0}^{[|(\ell-1)/2|]}(2^{2j}+2^{2(\ell-j-1)}-3)\\& &\times\sym\{(\vt_3(z)^{n_1+\cdots+n_{2j+1}}-\vt_2(z)^{n_1+\cdots+n_{2j+1}}-\vt_3(z)^{n_1+\cdots+n_{2j+1}})(\vt_3(z)^{n_{2j+2}+\cdots+n_{2\ell}}\\& &-\vt_2(z)^{n_{2j+2}+\cdots+n_{2\ell}}-\vt_4(z)^{n_{2j+2}+\cdots+n_{2\ell}})\}+[1-\frac{1}{6}\sum_{1\leq j_1, j_2;j_1+j_2\leq\ell-1}\left(\!\!\!\begin{array}{c}2\ell\\ 2j_1,2j_2\end{array}\!\!\!\right)\\& &+\frac{1}{2}\sum_{0\leq j_1,j_2;j_1+j_2\leq\ell-2}\left(\!\!\!\begin{array}{c}2\ell\\ 2j_1+1,2j_2+1\end{array}\!\!\!\right)-\frac{3}{2}\sum_{j=1}^{\ell-1}\left(\!\!\!\begin{array}{c}2\ell\\ 2j\end{array}\!\!\!\right)-\frac{1}{2}\sum_{i=0}^{\ell-1}\left(\!\!\!\begin{array}{c}2\ell\\ 2j+1\end{array}\!\!\!\right)(2^{2j}\\& &+2^{2(\ell-j-1)}-3)](\vt_2(z)^{n_1+\cdots+n_{2\ell}}+\vt_3(z)^{n_1+\cdots+n_{2\ell}}+\vt_4(z)^{n_1+\cdots+n_{2\ell}})
\\&=&\sum_{1\leq j_1\leq j_2\leq\ell-j_1-j_2}\sym\{h_{m_1+\cdots+m_{2j_1}+j_1}(z)h_{m_{2j_1+1}+\cdots+m_{2(j_1+j_2)}+j_2}(z)\\& &\times h_{m_{2(j_1+j_2)+1}+\cdots+m_{2\ell}+\ell-j_1-j_2}(z)\}
-2^8\Dlt_{24}(z)\sum_{1\leq j_1,j_2;j_1+j_2\leq\ell-2}\sym\{\rho_{m_1+\cdots+m_{2j_1+1}+j_1-1}(z)\\& &\times\rho_{m_{2j_1+2}+\cdots+m_{2(j_1+j_2+1)}+j_2-1}(z)h_{m_{2(j_1+j_2)+3}+\cdots+m_{2\ell}+\ell-j_1-j_2-1}(z)\}\\& &+3\sum_{j=1}^{[|\ell/2|]}\sym\{h_{m_1+\cdots+m_{2j}+j}(z)
 h_{m_{2j+1}+\cdots+m_{2\ell}+\ell-j}(z)\}+2^8\Dlt_{24}(z)\\& &\times\sum_{j=0}^{[|(\ell-1)/2|]}(2^{2j}+2^{2(\ell-j-1)}-3)\sym\{\rho_{m_1+\cdots+m_{2j+1}+j-1}(z)\rho_{m_{2j+2}+\cdots+m_{2\ell}+\ell-j-2}(z)\}\\& &+\left(4-\frac{2}{3}\sum_{j_1+j_2\leq\ell}\left(\!\!\!\begin{array}{c}2\ell\\ 2j_1,2j_2\end{array}\!\!\!\right)\right)h_{m_1+\cdots+m_{2\ell}+\ell}(z).\qquad\Box\hspace{4.7cm}(3.33)\end{eqnarray*}
\pse

According to the above theorem, we have the following  theta series of the lattice (3.9) with the assumption on $n_i$ as in the theorem:
$$\Theta_L(z)={1\over 2}h_{m_1}(z)\eqno(3.34)$$
when $k=1$;
$$\Theta_L(z)={1\over 2}h_{m_1+m_2+1}(z)-64\Dlt_{24}(z)\rho_{m_1-1}(z)\rho_{m_2-1}(z)\eqno(3.35)$$
when $k=2$;
\begin{eqnarray*}\Theta_L(z)&=&-{1\over 2}h_{m_1+m_2+m_3}(z)+{1\over 4}(h_{m_1}(z)h_{m_2+m_3}(z)+h_{m_2}(z)h_{m_1+m_3}(z)\\& &+h_{m_3}(z)h_{m_1+m_2}(z))-{1\over 8}h_{m_1}(z)h_{m_2}(z)h_{m_3}(z)\hspace{4.9cm}(3.36)\end{eqnarray*}
when $k=3$;
\begin{eqnarray*}&&\Theta_L(z)\\&= &16\Dlt_{24}(z)[2\:\sym\{\rho_{m_1-1}(z)\rho_{m_2+m_3+m_4}(z)\}-\sym\{h_{m_1+m_2+1}(z)\rho_{m_3-1}(z)\rho_{m_4-1}(z)\}]\\& &+\frac{1}{16}[3\:\sym\{h_{m_1+m_2+1}(z)h_{m_3+m_4+1}(z)\}-10h_{m_1+m_2+m_3+m_4+2}(z)]\hspace{2.5cm}(3.37)\end{eqnarray*}
when $k=4$;
\begin{eqnarray*}\Theta_L(z)&= &{1\over 32}[\mbox{sym}\{h_{m_1}(z)h_{m_2+m_3}(z)h_{m_4+m_5}(z)\}-\sym\{h_{m_1}(z)h_{m_2}(z)h_{m_3+m_4+m_5}(z)\}\\& &2\:\sym\{h_{m_1}(z)h_{m_2+m_3+m_4+m_5}(z)\}- \sym\{h_{m_1+m_2}(z)h_{m_3+m_4+m_5}(z)\}\\& &-4h_{m_1+m_2+m_4+m_5}(z)]\hspace{9.1cm}(3.38)\end{eqnarray*}
when $k=5$;
\begin{eqnarray*}& &\Theta_L(z)\\&= &{1\over 64}[-118h_{m_1+m_2+m_3+m_4+m_5+m_6+3}(z)+3\:\sym\{h_{m_1+m_2+1}(z)h_{m_3+m_4+m_5+m_6+2}(z)\}\\& &+\sym\{h_{m_1+m_2+1}(z)h_{m_3+m_4+1}(z)h_{m_5+m_6+1}(z)\}]-4\Dlt_{24}(z)[\sym\{h_{m_1+m_2+1}(z)\\& &\times \rho_{m_3-1}(z)\rho_{m_4+m_5+m_6}(z)\}+\sym\{h_{m_1+m_2+m_3+m_4+2}(z)\rho_{m_5-1}(z)\rho_{m_6-1}(z)\}\\& &-14\:\sym\{\rho_{m_1-1}(z)\rho_{m_2+m_3+m_4+m_5+m_6+1}(z)\}\\& &-5\:\sym\{\rho_{m_1+m_2+m_3}(z)\rho_{m_4+m_5+m_6}(z)\}]\hspace{7cm}(3.39)\end{eqnarray*}
when $k=6$. In particular, we have 
$$\mbox{the theta series of Niemeier lattice of type}\;D_{24}=E_4^3(z)+384\Dlt_{24}(z)\eqno(3.40)$$
by (3.34);
$$\mbox{the theta series of Niemeier lattice of type}\;D^2_{12}=E_4^3(z)-192\Dlt_{24}(z)\eqno(3.41)$$
by (3.35);
$$\mbox{the theta series of Niemeier lattice of type}\;D_8^3=E_4^3(z)-384\Dlt_{24}(z)\eqno(3.42)$$
by (3.36). We refer to [N] for the Niemeier lattices. 
\pse

{\bf Remark 3.2}. In the case $k=2\ell$ and $n_1=n_2=\cdots =n_{2\ell}=4$, the lattice $L$ in Theorem 3.1 contains a sublattice isomorphic to $R_{D_{8\ell}}$, and $\Theta_L(z)=h_{\ell}(z)/2$. The theta series $\Theta_L(z)$ in (3.12) is different from that in (3.11) if not all $m_i$ are zero.
\pse

We go back to the settings in (3.1)-(3.8). Assume that $k=4$ and $n_i=8m_i+4\es+2$ with $m_i\in\Bbb{N}$ and $\es=0,1$. Set
$$L=\Bbb{Z}(x^1_2+x^3_3+x^4_3)+ \Bbb{Z}(x^2_2+x^4_3+x^1_3)+\Bbb{Z}(x^3_2+x^1_3+x^2_3)+\Bbb{Z}(x^4_2+x^2_3+x^3_3)+Q\eqno(3.43)$$
(cf. (3.4)). 
\psp

{\bf Theorem 3.3}. {\it The lattice} $L$ {\it in (3.43) is a positive definite even unimodular lattice with the theta series}
\begin{eqnarray*}\Theta_L(z)&=&{1\over 2}h_{m_1+m_2+m_3+m_4+2\es+1}(z)-32\Dlt_{24}(z)(\rho_{m_1+m_2+\es-1}(z)\rho_{m_3+m_4+\es-1}(z)\\& &+\rho_{m_1+m_3+\es-1}(z)\rho_{m_2+m_4+\es-1}(z)+\rho_{m_1+m_4+\es-1}(z)\rho_{m_2+m_3+\es-1}(z)).\hspace{1.5cm}(3.44)\end{eqnarray*}

{\it Proof}. $L$ is unimodular by (1.11) and (3.7)-(3.9). By (3.5)-(3.7) and calculating the theta series of the cosets in $L/Q$, we obtain
\begin{eqnarray*}& &2^4\Theta_L(z)\\&=&\vt_2(z)^{n_3+n_4}(\vt_3(z)^{n_1}-\vt_4(z)^{n_1})(\vt_3(z)^{n_2}+\vt_4(z)^{n_2})
+\vt_2(z)^{n_1+n_4}(\vt_3(z)^{n_2}-\vt_4(z)^{n_2})\\
& &\times(\vt_3(z)^{n_3}+\vt_4(z)^{n_3})+\vt_2(z)^{n_1+n_2}(\vt_3(z)^{n_3}-\vt_4(z)^{n_3})(\vt_3(z)^{n_4}+\vt_4(z)^{n_4})\\& &+\vt_2(z)^{n_2+n_3}(\vt_3(z)^{n_4}-\vt_4(z)^{n_4})(\vt_3(z)^{n_1}+\vt_4(z)^{n_1})+\vt_2(z)^{n_1+n_3}(\vt_3(z)^{n_2}-\vt_4(z)^{n_2})\\&&\times(\vt_3(z)^{n_4}+\vt_4(z)^{n_4})+\vt_2(z)^{n_1+n_2+n_3+n_4}+\vt_2(z)^{n_2+n_4}(\vt_3(z)^{n_1}-\vt_4(z)^{n_1})\\& &\times(\vt_3(z)^{n_3}+\vt_4(z)^{n_3})+\vt_2(z)^{n_2+n_4}(\vt_3(z)^{n_3}-\vt_4(z)^{n_3})(\vt_3(z)^{n_1}+\vt_4(z)^{n_1})\\& &+\vt_2(z)^{n_1+n_2+n_3+n_4}+\vt_2(z)^{n_1+n_3}(\vt_3(z)^{n_4}-\vt_4(z)^{n_4})(\vt_3(z)^{n_2}+\vt_4(z)^{n_2})+\vt_2(z)^{n_2+n_3}\\& &\times(\vt_3(z)^{n_1}-\vt_4(z)^{n_1})(\vt_3(z)^{n_4}+\vt_4(z)^{n_4})+\vt_2(z)^{n_1+n_2}(\vt_3(z)^{n_4}-\vt_4(z)^{n_4})\\& &\times(\vt_3(z)^{n_3}+\vt_4(z)^{n_3})+\vt_2(z)^{n_1+n_4}(\vt_3(z)^{n_3}-\vt_4(z)^{n_3})(\vt_3(z)^{n_2}+\vt_4(z)^{n_2})\\& &+\vt_2(z)^{n_3+n_4}(\vt_3(z)^{n_2}-\vt_4(z)^{n_2})(\vt_3(z)^{n_1}+\vt_4(z)^{n_1})+(\vt_3(z)^{n_1}-\vt_4(z)^{n_1})(\vt_3(z)^{n_2}\\& &-\vt_4(z)^{n_2})(\vt_3(z)^{n_3}-\vt_4(z)^{n_3})(\vt_4(z)^{n_1}-\vt_4(z)^{n_4})+(\vt_3(z)^{n_1}+\vt_4(z)^{n_1})\\& &\times(\vt_3(z)^{n_2}+\vt_4(z)^{n_2})(\vt_3(z)^{n_3}+\vt_4(z)^{n_3})(\vt_4(z)^{n_1}+\vt_4(z)^{n_4})\\&=&\vt_2(z)^{n_3+n_4}[(\vt_3(z)^{n_1}-\vt_4(z)^{n_1})(\vt_3(z)^{n_2}+\vt_4(z)^{n_2})+(\vt_3(z)^{n_2}-\vt_4(z)^{n_2})\\& &\times(\vt_3(z)^{n_1}+\vt_4(z)^{n_1})]+\vt_2(z)^{n_1+n_4}[(\vt_3(z)^{n_2}-\vt_4(z)^{n_2})(\vt_3(z)^{n_2}+\vt_4(z)^{n_3})\\& &+(\vt_3(z)^{n_3}-\vt_4(z)^{n_3})(\vt_3(z)^{n_2}+\vt_4(z)^{n_2})]
+\vt_2(z)^{n_1+n_2}[(\vt_3(z)^{n_3}-\vt_4(z)^{n_3})\hspace{3cm}\end{eqnarray*}
\begin{eqnarray*}& &\times(\vt_3(z)^{n_4}+\vt_4(z)^{n_4})+(\vt_3(z)^{n_4}-\vt_4(z)^{n_4})(\vt_3(z)^{n_3}+\vt_4(z)^{n_3})]\\&&+\vt_2(z)^{n_2+n_3}[(\vt_3(z)^{n_1}-\vt_4(z)^{n_1})(\vt_3(z)^{n_4}+\vt_4(z)^{n_4})+(\vt_3(z)^{n_4}-\vt_4(z)^{n_4})\\& &\times(\vt_3(z)^{n_1}+\vt_4(z)^{n_1})]
+\vt_2(z)^{n_1+n_3}[(\vt_3(z)^{n_2}-\vt_4(z)^{n_2})(\vt_3(z)^{n_4}+\vt_4(z)^{n_4})\\& &+(\vt_3(z)^{n_4}-\vt_4(z)^{n_4})(\vt_3(z)^{n_2}+\vt_4(z)^{n_2})]+\vt_2(z)^{n_2+n_4}[(\vt_3(z)^{n_1}-\vt_4(z)^{n_1})\\& &\times(\vt_3(z)^{n_3}+\vt_4(z)^{n_3})+(\vt_3(z)^{n_3}-\vt_4(z)^{n_3})(\vt_3(z)^{n_1}+\vt_4(z)^{n_1})]
\\& &+2\vt_2(z)^{n_1+n_2+n_3+n_4}+(\vt_3(z)^{n_1}-\vt_4(z)^{n_1})(\vt_3(z)^{n_2}-\vt_4(z)^{n_2})(\vt_3(z)^{n_3}-\vt_4(z)^{n_3})\\& &\times(\vt_4(z)^{n_1}-\vt_4(z)^{n_4})+(\vt_3(z)^{n_1}+\vt_4(z)^{n_1})(\vt_3(z)^{n_2}+\vt_4(z)^{n_2})\\& &\times(\vt_3(z)^{n_3}+\vt_4(z)^{n_3})(\vt_4(z)^{n_1}+\vt_4(z)^{n_4})
\\&=&2[-(\vt_3(z)^{n_1+n_2}-\vt_2(z)^{n_1+n_2}-\vt_4(z)^{n_1+n_2})(\vt_3(z)^{n_3+n_4}-\vt_2(z)^{n_3+n_4}-\vt_4(z)^{n_3+n_4})\\& &-
(\vt_3(z)^{n_1+n_3}-\vt_2(z)^{n_1+n_3}-\vt_4(z)^{n_1+n_3})(\vt_3(z)^{n_2+n_4}-\vt_2(z)^{n_2+n_4}-\vt_4(z)^{n_2+n_4})
\\& &-(\vt_3(z)^{n_1+n_4}-\vt_2(z)^{n_1+n_4}-\vt_4(z)^{n_1+n_4})(\vt_3(z)^{n_2+n_3}-\vt_2(z)^{n_2+n_3}-\vt_4(z)^{n_2+n_3})\\& &
+4(\vt_2(z)^{n_1+n_2+n_3+n_4}+\vt_3(z)^{n_1+n_2+n_3+n_4}+\vt_4(z)^{n_1+n_2+n_3+n_4})]\\&=&2[4h_{m_1+m_2+m_3+m_4+2\es+1}-2^8\Dlt_{24}(z)\sym\{\rho_{m_1+m_2+\es-1}\rho_{m_3+m_4+\es-1}\}],\hspace{2.3cm}(3.45)\end{eqnarray*}
which implies (3.44).$\qquad\Box$
\psp

By (3.44) with $m_1=m_2=m_3=m_4=0$ and $\es=1$, we have 
$$\mbox{the theta series of Niemeier lattice of type}\;D^4_6=E_4(z)^3-480\Dlt_{24}(z)\eqno(3.46)$$
(cf. [N]).
\psp

{\bf Remark 3.4}. (a) The theta series in (3.44) is different from those in (3.11) and (3.12) if
$(m_1,m_2,m_3,m_4,\es)\neq (0,0,0,0,0)$. 

(b) The above families of lattices are the infinite families of positive definite even unimodular lattices whose theta series are  weighted symmetric polynomials of the functions of $\{\Dlt_{24}(z),h_n(z),\rho_n(z)\mid n\in\Bbb{Z}_+\}$ that we can find so far. We speculate that the theta series of the other infinite families of positive definite even-unimodular lattices may be related to the invariants of the other finite groups.

\vspace{0.5cm}

\noindent{\Large \bf References}

\hspace{0.5cm}

\begin{description}
\item[{[B1]}]
R. E. Borcherds, The Leech lattice and other lattices, Ph.D.
Dissertation, Univ. of Cambridge, 1983.

\item[{[B2]}]
---, Vertex algebras, Kac-Moody algebras, and the Monster,
{\it Proc. Natl. Acad. Sci. USA} {\bf 83} (1986), 3068-3071.

\item[{[C]}] ---, Three lectures on exceptional groups, in: {\it Finite
simple groups,} ed. by G. Higman and M. B. Powell, Chap.7, Academic
Pless, London/New York, 215-247, 1971.

\item[{[CS1]}] J. H. Conway and N. J. A. Sloane, Twenty-three
constructions of the Leech lattice, {\it Proc. Roy. Soc. London A} {\bf 
381} (1982), 275-283.

\item[{[CS2]}] ---, On the enumeration of
lattices of determinant one, {\it J. Number Theory} {\bf 15} (1982), 83-93.

\item[{[CS3]}] ---, {\it Sphere
Packings, Lattices and Groups,} Springer-Verlag, 1988.

\item[{[H]}] E. Hecke, Analytische arithmetik der positiven quadratischen formen, {\it Kgl. Danske Vid. Selskab. Mat.-fys. Medd.} {\bf 13} No. 12 (1940). 

\item[{[F]}] M. H. Freedman, The topolgy of 4-dimensional manifolds, {\it J. Diff. Geom.} {\bf 17} (1982), 357-453.

\item[{[FLM]}]
I. B. Frenkel, J. Lepowsky and A. Meurman, {\it Vertex Operator
Algebras and the Monster}, Pure and Applied Mathematics {\bf 134}, Academic Press Inc., Boston, 1988.

\item[{[L]}] D. H. Lehmer, Ramanujan's function $\tau(n)$, {\it Duke Math. J.} {\bf 10} (1943), 483-492.

\item[{[N]}] H. V. Niemeier, Definite quadratische formen der dimension 24 und diskriminante 1, {\it J. Number Theory} {\bf 5} (1973), 142-178.

\item[{[S1]}]
N. J. A. Sloane, codes over $GF(4)$ and complex lattices, {\it J. Algebra} {\bf 52} (1978), 168-181.

\item[{[S2]}]
---, Self-dual codes and lattices,
 {\it Proc. Symposia in Pure Math.} {\bf 34} (1979), 273-307.
  
\item[{[X1]}] X. Xu, Untwisted and twisted gluing techniquaes for constructing self-dual lattices, Ph.D. Dissertation, Rutgers University, 1992.

\item[{[X2]}]---, Self-dual lattices of type A, {\it Acta Math.} {\bf 175} (1995), 123-150.

\item[{[X3]}] ---, {\it Introduction to Vertex Operator Superalgebras and Their Modules}, Kluwer Academic Publishers, Dordrecht/Boston/London, 1998.

\end{description}
\end{document}